# Uniform stability estimates for the discrete Calderón problems

**S Ervedoza**[1,2] **and F de Gournay**[3]

1 CNRS ; Institut de Mathématiques de Toulouse UMR 5219 ; F-31062 Toulouse, France
2 Université de Toulouse ; UPS, INSA, INP, ISAE, UT1, UTM ; IMT ; F-31062 Toulouse, France,

E-mail: `ervedoza@math.univ-toulouse.fr`

3 Laboratoire de Mathématiques de Versailles, Université de Versailles Saint-Quentin-en-Yvelines, 45, avenue des Etats Unis, F-78035 Versailles Cedex, France.

E-mail: `gournay@math.uvsq.fr`

**Abstract.** In this article, we focus on the analysis of discrete versions of the Calderón problem in dimension $d \geq 3$. In particular, our goal is to obtain stability estimates for the discrete Calderón problems that hold uniformly with respect to the discretization parameter. Our approach mimics the one in the continuous setting. Namely, we shall prove discrete Carleman estimates for the discrete Laplace operator. A main difference with the continuous ones is that there, the Carleman parameters cannot be taken arbitrarily large, but should be smaller than some frequency scale depending on the mesh size. Following the by-now classical Complex Geometric Optics (CGO) approach, we can thus derive discrete CGO solutions, but with limited range of parameters. As in the continuous case, we then use these solutions to obtain uniform stability estimates for the discrete Calderón problems.

# 1. Introduction

The Calderón problem, or Electrical Impedance Tomography, amounts to retrieving the potential $q$ and the conductivity $\sigma$ from the knowledge of the current-to-voltage map on the boundary on a domain.

To be more precise, let $\Omega$ be a smooth bounded domain of $\mathbb{R}^d$. Given a conductivity $\sigma$ and a potential $q$, we define the operator $\Lambda[\sigma, q] : H^{1/2}(\partial\Omega) \to H^{-1/2}(\partial\Omega)$, the so-called Dirichlet-to-Neumann map, by

$$\Lambda[\sigma, q](g) = \sigma \nabla u[\sigma, q] \cdot n_{|\partial\Omega},$$

where $u[\sigma, q]$ solves the elliptic problem

$$\text{div}(\sigma \nabla u) + qu = 0 \text{ in } \Omega \text{ and } u = g \text{ on } \partial\Omega. \tag{1.1}$$

The Calderón problem is then the following: Given $\Lambda[\sigma, q]$, can we find the conductivity $\sigma$ and the potential $q$? Of course, the first question to be solved is the so-called uniqueness: is the map $(\sigma, q) \mapsto \Lambda[\sigma, q]$ injective?

Note that here, the knowledge of the Dirichlet-to-Neumann map amounts to know the answer of the system (1.1) to *any* boundary data $g \in H^{1/2}(\partial\Omega)$.

Before recalling the known results on the above Calderón problem, let us point out the by-now well known fact that when $\sigma$ is a scalar, the Liouville transform allows us to rewrite the Calderón problem as follows: Given the map $\Lambda[1, q] : H^{1/2}(\partial\Omega) \to H^{-1/2}(\partial\Omega)$, that is the map defined by

$$\Lambda[1, q](g) = \partial_n u[q]_{|\partial\Omega},$$

where $u[q]$ is the solution of

$$\Delta u + qu = 0 \text{ in } \Omega \text{ and } u = g \text{ on } \partial\Omega, \tag{1.2}$$

can we reconstruct $q$?

In the sequel, we shall focus on this very precise Calderón problem, based on the Laplace operator in (1.2). Therefore, from now on, we shall denote $\Lambda[1, q]$ by $\Lambda[q]$ without any confusion.

Of course, here again, the uniqueness is the first question to ask, that is the study of the injectivity of the map $\Lambda : q \mapsto \Lambda[q]$, and then the stability of the inverse, namely trying to understand the modulus of continuity of the map $\Lambda^{-1}$.

There is an extensive literature on the Calderón problem. First, some energy considerations solve the problem of uniqueness of $\sigma$ given $\Lambda[\sigma, 0]$ if $\sigma$ is supposed to be scalar (isotropic) and piece-wise analytic, see [19].

For isotropic conductivities, a great deal of work has been spent in order to reduce the regularity hypothesis on $\sigma$ in order to guarantee the injectivity of the mapping $\sigma \mapsto \Lambda[\sigma, 0]$. In dimension 2, using techniques of complex analysis, the conductivity is needed to be $L^\infty(\Omega)$ only (see [2, 10]). In dimension $d \geq 3$, using the CGO technique that will be mimicked hereafter, the conductivities must lie in the Sobolev space $W^{3/2,p}(\Omega)$ for $p > 2d$ (see [22, 9]) in order to ensure uniqueness. This effect is not totally surprising, since unique continuation (which seems closely related) behaves differently in dimension 2 and in greater dimension.

One may also suppose that the operator $\Lambda[q]$ is known only on a part of the boundary. In this case uniqueness and stability still hold under some suitable "illuminations" conditions if $q \in L^\infty(\Omega)$ (or $\sigma$ in $W^{2,\infty}(\Omega)$ by the Liouville transform) using the so-called limiting Carleman weights (see [18, 14]).

Concerning anisotropic conductivities, the celebrated Tartar's counterexample show that two conductivities differing only by the pull-back of a diffeomorphism equal to the identity on the boundary yield the same Dirichlet-to-Neumann map. In dimension 2 and anisotropic $\sigma$, Tartar's counterexample has been proven to be the only obstruction to uniqueness using complex analysis in isothermal coordinates (see [21, 3]) for $L^\infty(\Omega)$ conductivities. The case of the dimension $d \geq 3$ is not yet fully understood, the case of real analytic anisotropic conductivities is solved in [20] and encouraging progress has been done recently if the conductivity exhibits a preferred direction that allows to treat the geometry like a transformation of a cylinder (see [13]). To sum up the anisotropic case, the 2-d case is rather well understood, mainly thanks to complex analysis, whereas higher dimensions are more of a *Terra Incognita*.

The question we are interested in is the following. Discretizing the Calderón problem, can we still get stability estimates uniformly with respect to the mesh size?

To be more precise, we consider a discrete version of the elliptic equation (1.2). Let $\mathcal{W}_h$ be a mesh that approximates the set $\Omega$ and, for $g_h \in H_h^{1/2}(\partial\mathcal{W}_h)$, let $u_h[q_h]$ be the solution of

$$\Delta_h u_h + q_h u_h = 0 \text{ in } \mathring{\mathcal{W}}_h \text{ and } u_h = g_h \text{ on } \partial\mathcal{W}_h. \tag{1.3}$$

Here, the index $h$ means that we are solving a discrete approximation of (1.2) with parameter $h > 0$, corresponding for instance to the mesh size. Hence, $\mathring{\mathcal{W}}_h$ is an approximation of $\Omega$ (we shall come back later on the meaning of $\mathring{\mathcal{W}}_h$ compared to $\mathcal{W}_h$, see (2.1) for more details), $\Delta_h$ is an approximation of the Laplace operator $\Delta$, and the space $H_h^{1/2}(\partial\mathcal{W}_h)$ is the discrete counterpart of $H^{1/2}(\partial\Omega)$. We shall not give more precise descriptions of these approximations at this stage, and we refer the reader to the rest of the article for more details.

We can then introduce the discrete DtN maps $\Lambda_h$ defined by

$$\Lambda_h[q_h](g_h) = \partial_{n,h} u_h[q_h]_{|\partial\mathcal{W}_h},$$

where $u_h[q_h]$ is the solution of (1.3) corresponding to the Dirichlet data $g_h$ and $\partial_{n,h} u_h[q_h]$ is an approximation of the normal derivative of $u_h[q_h]$, see (2.15) for an explicit description of that quantity.

One should expect that, if the discrete approximations (1.3) converge (as $h \to 0$) to the continuous equation (1.2), the maps $q_h \mapsto \Lambda_h[q_h]$ as well as their reconstruction $\Lambda_h[q_h] \mapsto q_h$ should also converge, in some sense, to their continuous counterpart $q \mapsto \Lambda[q]$ and $\Lambda[q] \mapsto q$. Of course, these naive statements need to be thoroughly study, as it turns out that the maps $(\Lambda_h)_{h>0}$ all are non-linear, and therefore their inverses may have very different behavior. Our goal is to make a first attempt at describing the discrete DtN maps and their stability results in the case $d \geq 3$. In order to obtain convergence of the discrete reconstruction towards the continuous one, we think that these estimates, that will be derived uniformly with respect to $h$, are a step of crucial importance, see our comments in Section 7.

One of the difficulty to get uniform stability estimates is that the continuous approach, based on CGO, strongly makes use of arbitrary large frequency solutions, whereas in a discrete case of mesh-size $h > 0$, the frequency of solutions is limited to the scale $1/h$.

Another restriction linked to the discretization of the equation is that the convergence of the spectrum is only guaranteed for frequencies smaller than $h^{-2/3}$.

For instance, the discrete Laplacian in 1-d on a uniform mesh,

$$(\Delta_h u_h)_j = \frac{1}{h^2} \left(u_{j+1,h} - 2u_{j,h} + u_{j-1,h}\right), \quad j \in \{1, \cdots, 1/h - 1\}, \quad u_{0,h} = u_{1/h,h} = 0,$$

admits as eigenvalues,

$$\lambda_{k,h} = \frac{4}{h^2} \sin^2\left(\frac{k\pi h}{2}\right), \quad k \in \{1, \cdots, 1/h - 1\},$$

whereas the eigenvalues of the Laplacian with Dirichlet boundary conditions on $(0,1)$ are given by $\lambda_k = k^2 \pi^2$. Hence

$$\left|\sqrt{\lambda_{k,h}} - \sqrt{\lambda_k}\right| \leq C k^3 h^2, \quad \forall k \in \{1, \cdots, 1/h - 1\},$$

meaning that the convergence of the spectrum is achieved only in the range $k = o(h^{-2/3})$, corresponding to $\sqrt{\lambda_{k,h}} = o(h^{-2/3})$ or frequencies of order $o(h^{-2/3})$.

Since the CGO solutions are by their nature, linked to eigenvalues of the Laplacian, when speaking convergence of the discrete Calderón problems towards the continuous one, these scales of frequency should be kept in mind. Above these scales, one can expect very surprising behavior between the continuous setting and the discrete ones.

Let us point out a third and complex effect of the discretization, which is glaring in the exact controllability problem for the wave equation, where the discrete controls of smallest norm do not converge to the continuous one, see [16, 24].

This effect is due to numerical dispersion that introduces high-frequency spurious waves of null speed, despite the fact every solution of the continuous wave equation travels at velocity one (see [23]). Since these spurious waves do not propagate in the domain, they cannot be controlled. We refer the reader interested by these questions to the seminal articles [15, 17] and the subsequent survey article [24].

Of course, CGO solutions correspond to plane waves and therefore discrete versions of CGO solutions may have very different behaviors as in the continuous case.

Our last remark is based on counting the unknowns and the equations. For each $h > 0$, the discrete Calderón problem is stated in a finite-dimensional space: If the number of points in a direction is equivalent to $N$, knowing the whole DtN map is tantamount to knowing $N^{d-1}$ potentials, each of which yielding a current, that is $N^{d-1}$ data, on the boundary. Hence the knowledge of the DtN map is equivalent to the knowledge of $N^{2d-2}$ data, whereas the number of unknown is proportional to $N^d$. Therefore, in dimension $d \geq 3$, there are way too many data for the number of unknown and the problem is ill-posed in the sense that they might not exist a solution at all for an arbitrary set of data. Even if uniqueness of a solution is guaranteed, existence is here of the essence. This remark seems to indicate that solving the discrete Calderón problem is an easy matter, but we emphasize that the map $\Lambda_h$ is non-linear and lies in a finite-dimensional space of very large dimension, going to infinity as $h \to 0$, thus making the resolution of the discrete Calderón problem a challenging issue.

The results in the continuous case and in the discrete case are in a contrast so sharp that it is on the verge of being paradoxical. In the discrete case, iterative reconstruction algorithms exists (see the review [5]), based on continuous analysis and a direct algorithm is proposed in [12] which strongly uses the intrinsic finite-dimensional nature of these systems. However, these latest reconstruction algorithms use highly oscillating test functions. As a consequence, for a set of data that represents a DtN map (possibly with some errors), the reconstructed conductivities and potentials are designed to match the $N^d$ data of higher frequency. Let us also point out the direct algorithm of [4] that builds the set of point on which the conductivity is reconstructed accordingly to the data, but, as far as we know, this reconstruction algorithm does not yield any uniform stability estimates for the discrete Calderón problem, where uniform stands for uniformly with respect to the discretization parameter $h > 0$.

To study the discrete Calderón problems, we shall then develop similar tools as for the continuous Calderón problem in dimension $d \geq 3$. To be more precise, our analysis will be based on the construction of discrete CGO solutions, solution of (1.3) which are close to the usual harmonic functions. Of course, this cannot be done for any arbitrarily large frequencies but rather for a limited range of frequency, depending on the mesh-size $h$.

One of the main difficulties we shall encounter to construct these solutions is to develop a discrete Carleman estimate for the elliptic problem (1.3). Of course, here again, as expected in view of the above discussion, the parameters in the Carleman estimate should be limited in some range depending on the mesh size (roughly speaking, the parameter in the Carleman estimate correspond to the frequency of the CGO solutions).

Note that, in this step, the recent works [6, 7, 8] are worth mentioning: to our knowledge, they are so far the only works in which discrete Carleman estimates have been derived, but having in mind the study of the exact controllability of the discrete heat equations.

In Section 2, we define the discrete problems, introduce some notations and state our main result. In Section 3, we prove a Carleman estimate which is the basis of the proof of existence of the CGO solutions, given in Section 4. In Section 5, we prove the stability estimate in terms of the different scales that come into play in the discrete Calderón problems. In Section 6, we focus on the case of uniform meshes, in which stronger stability estimates can be proved. In Section 7, we give some further comments and open problems.

## 2. Definition of the discrete problem

In this section, we introduce notations specific to the discrete problem that will be used throughout this paper.

Set $\Omega \subset [0,1]^d$, $d \geq 3$, a domain meshed with $\mathcal{M}_h$, a finite subset of $\bar{\Omega}$. For $N \in \mathbb{N}$, define $\mathcal{K}_h$ the regular cartesian grid of $[0,1]^d$:

$$\mathcal{K}_h := \left\{ x \in [0,1]^d \text{ such that } \exists k \in [\![0, N-1]\!]^d \text{ such that } x = \frac{k}{N} \right\}.$$

Suppose that the mesh $\mathcal{M}_h$ is a perturbation of $\mathcal{W}_h$, a subset of $\mathcal{K}_h$, in the sense that there exists a smooth diffeomorphism $F : [0,1]^d \to [0,1]^d$ close to the identity such that $F(\mathcal{W}_h) = \mathcal{M}_h$.

When treating discrete problems where nodes of $\mathcal{M}_h$ are directly related to degrees of freedom of the linear system (Lagrangian finite elements or finite differences), one usually introduce the concept of "connection" between two nodes (or "neighbours") of $\mathcal{W}_h$. We use the notation $x \sim y$ to denote the fact that two points $x \neq y \in \mathcal{W}_h$ are neighbours. We shall suppose that the set of connections or edges is the same for every point and is symmetric, that is:

**Assumption 1** *Let $h > 0$ be the size of the mesh. There exists $k$ vectors $(e_i)_{i=1\ldots k} \in \mathbb{Z}^d$ such that $e_i \neq -e_j$ for all $i \neq j$ and such that for all $x, y$ in $\mathcal{W}_h$,*

$$x \sim y \iff \exists i \text{ such that } x - y = he_i \text{ or } y - x = he_i.$$

*We further assume that the set of vectors $(e_i)_{i=1\ldots k}$ spans $\mathbb{R}^d$.*

If $x = y + he_i$, $x$ is called the "neighbour of $y$ in direction $e_i$", whereas $y$ is the neighbour of $x$ in direction $-e_i$. Note that, since the set of connections is the same for every point, the mesh has to have some translation invariance that restricts drastically the choice of possible meshes to, namely, cartesian grids.

Also note that we did not assume $k = d$. In particular, $k$ can be much larger than $d$. In the case $k = d$, much more can be said, see Section 3.6 and Section 6.

*2.1. Dual meshes and operators*

For any set of points $\mathcal{A}_h$, define $\mathcal{A}_h^i$, the dual set in the direction $e_i$, as

$$\mathcal{A}_h^i := \left\{ x + \frac{h}{2} e_i \text{ such that } x \in \mathcal{A}_h \right\} \cap \left\{ x - \frac{h}{2} e_i \text{ such that } x \in \mathcal{A}_h \right\}.$$

As an immediate corollary, $y \in \mathcal{A}_h^i$ if and only if both $y + \frac{h}{2} e_i$ and $y - \frac{h}{2} e_i$ are in $\mathcal{A}_h$, hence $y \in \mathcal{A}_h^i$ is the middle point of a segment connecting two neighbours in direction $e_i$ of points of $\mathcal{A}_h$. Similarly, define the sets $\mathcal{A}_h^{ij} = (\mathcal{A}_h^i)^j = \mathcal{A}_h^{ji}$. The inclusion $\mathcal{A}_h^{ii} \subset \mathcal{A}_h$ is strict, since $x \in \mathcal{A}_h^{ii}$ if and only if $x \in \mathcal{A}_h$ and both its neighbours in direction $e_i$ and $e_{-i}$ are in $\mathcal{A}_h$. Define the interior and the boundary of a set $\mathcal{A}_h$ as

$$\mathring{\mathcal{A}}_h := \bigcap_{i=1}^k \mathcal{A}_h^{ii} \text{ and } \partial^i \mathcal{A}_h := \mathcal{A}_h \setminus \mathcal{A}_h^{ii} \text{ and } \partial \mathcal{A}_h := \mathcal{A}_h \setminus \mathring{\mathcal{A}}_h. \tag{2.1}$$

For any finite set of point $\mathcal{A}_h$, denote as $C(\mathcal{A}_h)$ the set of functions from $\mathcal{A}_h$ to $\mathbb{C}$ and define the average and difference operators in the direction $e_i$ as the operators from $C(\mathcal{A}_h)$ to $C(\mathcal{A}_h^i)$ :

$$a_i(u_h) := y \mapsto \frac{1}{2} \left( u_h(y + \frac{h}{2} e_i) + u_h(y - \frac{h}{2} e_i) \right),$$

$$d_i(u_h) := y \mapsto \frac{1}{h} \left( u_h(y + \frac{h}{2} e_i) - u_h(y - \frac{h}{2} e_i) \right).$$

These operators are the discrete version of the average and the directional derivation. Note that they depend on $h > 0$ and should rather be denoted by $a_{i,h}$ and $d_{i,h}$ respectively, but there, no confusion can occur and we remove it with a slight abuse of notations.

We shall also sometimes denote by $d$ generic difference operators which coincides with one of the $d_i$ and $a$ generic average operator which coincides with one of the $a_i$.

Here are a list of easy lemmas that will be used thereafter.

**Lemma 2.1** *As operators from $C(\mathcal{A}_h)$ to $C(\mathcal{A}_h^{ij})$, we have the following identities:*

$$a_i a_j = a_j a_i, \quad a_i d_j = d_j a_i, \quad d_i d_j = d_j d_i.$$

Moreover, for any $(u_h, v_h) \in C(\mathcal{A}_h)^2$,

$$a_i(u_h v_h) = a_i(u_h)a_i(v_h) + \frac{h^2}{4}d_i(u_h)d_i(v_h) \quad \text{on } \mathcal{A}_h^i, \tag{2.2}$$

$$d_i(u_h v_h) = d_i(u_h)a_i(v_h) + a_i(u_h)d_i(v_h) \quad \text{on } \mathcal{A}_h^i, \tag{2.3}$$

$$u_h = a_i a_i(u_h) - \frac{h^2}{4}d_i d_i(u_h) \quad \text{on } \mathcal{A}_h^{ii}. \tag{2.4}$$

**Lemma 2.2** *Denote the hyperbolic cosinus and sinus as* ch *and* sh, *respectively. For any $u_h \in C(\mathcal{A}_h)$ and $i \in \{1, \cdots, k\}$, on $\mathcal{A}_h^i$ we have*

$$a_i(e^{u_h}) = e^{a_i u_h}\text{ch}\left(\frac{h d_i u_h}{2}\right), d_i(e^{u_h}) = \frac{2}{h}e^{a_i u_h}\text{sh}\left(\frac{h d_i u_h}{2}\right),$$

$$a_i(\text{sh}(u_h)) = \text{sh}(a_i u_h)\text{ch}\left(\frac{h d_i u_h}{2}\right),$$

$$d_i(\text{sh}(u_h)) = \frac{2}{h}\text{ch}(a_i u_h)\text{sh}\left(\frac{h d_i u_h}{2}\right),$$

$$a_i(\text{ch}(u_h)) = \text{ch}(a_i u_h)\text{ch}\left(\frac{h d_i u_h}{2}\right),$$

$$d_i(\text{ch}(u_h)) = \frac{2}{h}\text{sh}(a_i u_h)\text{sh}\left(\frac{h d_i u_h}{2}\right).$$

*Besides, for any vector $s \in \mathbb{R}^d$ and $i \in \{1, \cdots, k\}$, $d_i(s \cdot x) = s \cdot e_i$.*

**Lemma 2.3** *For any $i \in \{1, \cdots, k\}$, $u_h \in C(\mathcal{A}_h)$, $z_h \in C(\mathcal{A}_h^i)$, we have*

$$d_i(z_h(a_i u_h)) = a_i(z_h(d_i u_h)) + (d_i z_h)u_h \text{ on } \mathring{\mathcal{A}}_h, \tag{2.5}$$

$$a_i(z_h(a_i u_h)) = (a_i z_h)u_h + \frac{h^2}{4}d_i(z_h(d_i u_h)) \text{ on } \mathring{\mathcal{A}}_h. \tag{2.6}$$

Lemmas 2.1–2.2–2.3 can all be derived easily by straightforward computations left to the reader.

## 2.2. The discrete Laplace operator

On $\mathcal{W}_h$, we are looking for an approximation of the solution $u \in H^1(\Omega)$ of the continuous problem (1.2).

The discrete problem then amounts to finding $u_h \in C(\mathcal{W}_h)$ solution of (1.3) for some discrete operator $\Delta_h$ that approximates $\Delta$ on $\mathcal{W}_h$. Of course, this operator $\Delta_h$ is defined on a finite-dimensional space and can then be considered as a matrix, as it is usually done when doing numerics. But one can also see this operator $\Delta_h$ as an operator from $C(\mathcal{W}_h)$ to $C(\mathring{\mathcal{W}}_h)$, and we shall adopt this point of view in the rest of the article.

We shall assume the following property:

**Assumption 2** *For all $h > 0$, there exists $k$ functions $(\sigma_h^i)_{i=1\ldots k} \in C(\mathcal{W}_h^i)$, such that*

$$\Delta_h u_h = \sum_{i=1}^k d_i(\sigma_h^i d_i u_h) \text{ on } \mathring{\mathcal{W}}_h, \quad \forall u_h \in C(\mathcal{W}_h).$$

For some standard problems, the existence of such $\sigma_h^i$ is made in appendix, it relies solely on the following properties:

- There is a one to one correspondence between the linear system and the nodes.
- The linear system is symmetric.
- The function constant equal to 1 is in the kernel of $\Delta_h$ on $\mathring{\mathcal{W}}_h$.

The value of $\sigma_h^i(x)$ is given in the appendix for some numerical methods. Roughly speaking, it corresponds to the weight that the numerical method gives to the connection of the nodes $(x + he_i, x - he_i)$. These functions $\sigma_h^i$ encode the information on the numerical method that has been chosen.

Given a set $\mathcal{A}_h$, for $u_h \in C(\mathcal{A}_h)$, we define its $L^\infty(\mathcal{A}_h)$ norm as

$$\|u_h\|_{L^\infty(\mathcal{A}_h)} := \max_{x \in \mathcal{A}_h}(|u_h(x)|).$$

We suppose that the mesh regularity and the choice of discretization yield the following properties on $\sigma_h^i$:

**Assumption 3** *Suppose that $\sigma_h^i$ can be extended to a function of $C(\mathcal{K}_h^i)$. Still denoting this extension by $\sigma_h^i$, we define:*

$$\epsilon_d(h) = \sum_{i,j} \|d_j(\sigma_h^i)\|_{L^\infty(\mathcal{K}_h^i)}, \tag{2.7}$$

$$\epsilon_a(h) = \|\sum_{i=1}^{k} a_i(\sigma_h^i) e_i \otimes e_i - Id\|_{L^\infty(\mathring{\mathcal{K}}_h)}, \tag{2.8}$$

$$M(h) = \sum_{i} \|\sigma_h^i\|_{L^\infty(\mathcal{K}_h^i)}, \tag{2.9}$$

*and we suppose that*

$$M := \sup_{h \to 0} M(h) < \infty.$$

We shall compute explicitly $\sigma$ for some standard examples in the appendix and show where these two scalings $\epsilon_a$ and $\epsilon_d$ come from. Heuristically, condition (2.8) is an "isotropic" condition, hence the notation $\epsilon_a$, where $a$ stands for "anisotropy". The scaling $\epsilon_d$ comes from the regularity of the diffeomorphism $F$ that modifies the domain $\mathcal{W}_h$ into $\mathcal{M}_h$.

Note that Assumption 3 extends the operator $\Delta_h$ on the whole mesh $\mathcal{K}_h$. Assumption 3 therefore concerns the regularity of $\mathcal{W}_h$. Indeed, extending $\sigma_h^i$ to a function of $C(\mathcal{K}_h^i)$ with the required properties implies that the domain $\Omega$, approximated by $\mathcal{W}_h$, is regular.

We shall also make and additional assumption on the geometry of the domain $\mathcal{W}_h$:

**Assumption 4** *Denote by $\ddot{\mathcal{K}}_h$ the interior of the interior of $\mathcal{K}_h$ (points that are at least two nodes away form $\partial \mathcal{K}_h$). We suppose that there exists a domain $\mathcal{B}_h$ such that $\mathcal{W}_h \subset \mathcal{B}_h$, $\mathcal{B}_h \subset \ddot{\mathcal{K}}_h$ and such that there exists $\psi \in C_c(\mathcal{B}_h)$ such that $\psi = 1$ on $\mathcal{W}_h$ and $\|d_i(\psi)\|_{L^\infty(\mathcal{K}_h)} \leq M_0$ for all $i \in \{1, \cdots, k\}$ and $\|\Delta_h \psi\|_{L^\infty(\mathcal{K}_h)} \leq M_0$, where $M_0$ is a constant independent of $h$.*

Note that $\psi$ in Assumption 4 is assumed to be independent of $h > 0$. This can indeed be done by assuming that $\bar{\Omega} \subset (0,1)^d$ and constructing a function $\psi$ compactly supported in $(0,1)^d$ and equal to one on $\bar{\Omega}$, and then taking $h > 0$ small enough.

Assumption 4 states that the domain $\mathcal{W}_h$ is at a fixed distance of the boundary $\partial \mathcal{K}_h$, in particular, we have

$$\mathcal{W}_h \subset \mathring{\mathcal{B}}_h \subset \mathcal{B}_h \subset \ddot{\mathcal{K}}_h \subset \mathcal{K}_h. \tag{2.10}$$

*2.3. Integrals and Green's identity*

For any set of point $\mathcal{A}_h$, for any $(u_h, v_h) \in C(\mathcal{A}_h)$, we define the following quantities:

$$\int_{\mathcal{A}_h} u_h(x) := h^d \sum_{x \in \mathcal{A}_h} u_h(x), \quad (u_h, v_h)_{\mathcal{A}_h} := \int_{\mathcal{A}_h} u_h(x) v_h(x). \quad (2.11)$$

Note that $(\cdot, \cdot)_{\mathcal{A}_h}$ defines a scalar product on the functions of $C(\mathcal{A}_h)$.

Remark that we choose the scaling $h^d$ of the volume integral. Therefore, in the sequel (and especially in Green's identity), the physical scaling of boundary integrals will be $h^{-1} \int_{\partial \mathcal{A}_h} u$.

We also define the following norms in $C(\mathcal{A}_h)$:

$$\|u_h\|_{L^2(\mathcal{A}_h)} := \sqrt{(u_h, u_h)_{\mathcal{A}_h}}, \quad \|u_h\|_{\dot{H}^1(\mathcal{A}_h)} := \left( \sum_i \int_{\mathcal{A}_h^i} \sigma_h^i d_i(u_h)^2 \right)^{1/2},$$

$$\|u_h\|_{H^1(\mathcal{A}_h)} := \|u_h\|_{L^2(\mathcal{A}_h)} + \|u_h\|_{\dot{H}^1(\mathcal{A}_h)}. \quad (2.12)$$

We are now in position to state the following discrete version of Green's identity:

**Proposition 2.4 (Integration by parts)** *Define $C_c(\mathcal{A}_h)$ the space of functions with compact support as*

$$C_c(\mathcal{A}_h) := \{u_h \in C(\mathcal{A}_h) \text{ s.t. } u_h = 0 \text{ on } \partial \mathcal{A}_h\}.$$

*For any $v_h \in C(\mathcal{A}_h^i)$, $u_h \in C_c(\mathcal{A}_h)$,*

$$\int_{\mathcal{A}_h^i} d_i(u_h) v_h = -\int_{\mathring{\mathcal{A}}_h} u_h d_i(v_h), \quad \int_{\mathcal{A}_h^i} a_i(u_h) v_h = \int_{\mathring{\mathcal{A}}_h} u_h a_i(v_h). \quad (2.13)$$

**Proof** Even if the proof of Green's identity is standard, we show it here in order to work out the notations. Below we only prove Green's identity for the difference operator and not for the average one, which can be proved similarly and is left to the reader. Let us introduce the set

$$\mathcal{A}_h^{\pm, i} = \left\{ x \text{ s.t. } x \pm \frac{h}{2} e_i \in \mathcal{A}_h^i \right\},$$

yielding in particular $\mathcal{A}_h^{ii} = \mathcal{A}_h^{+,i} \cap \mathcal{A}_h^{-,i}$. Moreover, performing two discrete change of variables,

$$\int_{\mathcal{A}_h^i} d_i(u_h) v_h = h^d \sum_i \sum_{x \in \mathcal{A}_h^i} \frac{1}{h} \left( u_h(x + \frac{h}{2} e_i) - u_h(x - \frac{h}{2} e_i) \right) v_h(x)$$

$$= h^d \sum_i \left( \sum_{y \in \mathcal{A}_h^{-,i}} \frac{1}{h} u_h(y) v_h(y - \frac{h}{2} e_i) - \sum_{y \in \mathcal{A}_h^{+,i}} \frac{1}{h} u_h(y) v_h(y + \frac{h}{2} e_i) \right).$$

Since $u_h = 0$ on $\partial \mathcal{A}_h$, all the sums on $\mathcal{A}_h^{\pm, i}$ are in fact on the set $\mathring{\mathcal{A}}_h \subset \mathcal{A}_h^{\pm, i}$, and hence

$$\int_{\mathcal{A}_h^i} d_i(u_h) v_h = h^d \sum_i \sum_{y \in \mathring{\mathcal{A}}_h} \frac{1}{h} u_h(y) \left( v_h(y - \frac{h}{2} e_i) - v_h(y + \frac{h}{2} e_i) \right)$$

$$= -\int_{\mathring{\mathcal{A}}_h} u_h d_i(v_h).$$

This concludes the proof of Proposition 2.4. □

In order to get Green's identity with boundary term, we have to define, for each $i \in \{1, \cdots k\}$, the exterior normal to the set $\mathcal{W}_h \subset \mathcal{K}_h$ in direction $e_i$ as $n_h^i \in C(\partial^i \mathcal{W}_h)$ :

$$\forall x \in \partial^i \mathcal{W}_h, \quad n_h^i(x) = \begin{cases} 1 & \text{if } x - \frac{h}{2}e_i \in \mathcal{W}_h^i \text{ but } x + \frac{h}{2}e_i \notin \mathcal{W}_h^i, \\ -1 & \text{if } x + \frac{h}{2}e_i \in \mathcal{W}_h^i \text{ but } x - \frac{h}{2}e_i \notin \mathcal{W}_h^i, \\ 0 & \text{if } x + \frac{h}{2}e_i \notin \mathcal{W}_h^i \text{ and } x - \frac{h}{2}e_i \notin \mathcal{W}_h^i. \end{cases}$$

This definition makes sense since $x + \frac{h}{2}e_i \in \mathcal{W}_h^i$ and $x - \frac{h}{2}e_i \in \mathcal{W}_h^i$ would imply that $x \in \mathcal{W}_h^{ii}$ hence $x \notin \partial^i \mathcal{W}_h$. Remark also that $n^i(x) = 0$ on $\partial^i \mathcal{W}_h$ if and only if $x$ does not have any neighbours in the direction $e_i$ nor in the direction $-e_i$, so that $x$ is an "isolated" point in the $\pm e_i$ direction.

**Proposition 2.5** *Denote by $I_h$ the operator from $C(\mathcal{W}_h^i)$ to $C(\mathcal{K}_h^i)$ that extends functions by zero outside $\mathcal{W}_h^i$, then*

$$d_i(I_h(v_h)) = -\frac{2}{h}a_i(I_h(v_h))n_h^i \quad \text{on } \partial^i \mathcal{W}_h. \tag{2.14}$$

*Extending $n_h^i$ by 0 outside $\partial^i \mathcal{W}_h$, for all $u_h \in C(\mathcal{W}_h)$, we define $\partial_{n,h} u_h \in C(\partial \mathcal{W}_h)$ as*

$$\partial_{n,h} u_h = \sum_{i=1}^k 2(a_i \circ I)(\sigma^i d_i u_h) n_h^i. \tag{2.15}$$

*With this definition, we have the following formula: for all $w_h \in C(\mathcal{W}_h)$,*

$$\frac{1}{h}\int_{\partial \mathcal{W}_h}(\partial_{n,h} u_h) w_h = \int_{\mathring{\mathcal{W}}_h}(\Delta_h u_h)w_h + \sum_{i=1}^k \int_{\mathcal{W}_h^i} \sigma^i(d_i u_h)(d_i w_h). \tag{2.16}$$

Let us emphasize that the factor $1/h$ in front of the integral in $\partial \mathcal{W}_h$ comes from the fact that this latter integral is defined as $\int_{\partial \mathcal{W}_h} u_h = h^d \sum_{x \in \partial \mathcal{W}_h} u_h(x)$ and not with the standard scaling for boundaries which is $h^{d-1}$.

**Proof** We first prove equation (2.14). Remark that if $\mathcal{W}_h \subset \mathring{\mathcal{K}}_h$, then $\mathcal{W}_h^i \subset \mathcal{K}_h^i$ so that the operator $I_h$ can indeed be defined. Moreover, the functions $d_i(I_h(v_h))$ and $a_i(I_h(v_h))$ both belong to $C(\mathcal{K}_h^{ii})$. Due to the inclusions

$$\partial^i \mathcal{W}_h \subset \mathcal{W}_h \subset \mathring{\mathcal{K}}_h \subset \mathcal{K}_h^{ii},$$

it makes sense to look at the values of $d_i(I_h(v_h))$ on $\partial^i \mathcal{W}_h$. Set $x \in \partial^i \mathcal{W}_h$, then either $n_h^i(x) = -1, 0$ or $1$. Since the three cases are treated the same way, we may suppose that $n_h^i(x) = 1$, that is $x - \frac{h}{2}e_i \in \mathcal{W}_h^i$ and $x + \frac{h}{2}e_i \notin \mathcal{W}_h^i$. In this case $(I_h v_h)(x + \frac{h}{2}e_i) = 0$ and

$$d_i(I_h v_h)(x) = \frac{1}{h}(-v_h(x - \frac{h}{2}e_i)) = -\frac{2}{h}a_i(I_h v_h)(x),$$

which is the sought result.

We now turn our attention to proving Green's formula by itself. Denote by $\tilde{I}_h$ the operator from $C(\mathcal{W}_h)$ to $C(\mathcal{K}_h)$ that extends functions by zero on $\mathcal{K}_h \setminus \mathcal{W}_h$. We have

$$d_i(\tilde{I}_h(w_h)) = d_i(w_h) \text{ on } \mathcal{W}_h^i \text{ and } I_h(\sigma_h^i d_i u_h) = 0 \text{ outside } \mathcal{W}_h^i.$$

Hence

$$\sum_{i=1}^k \int_{\mathcal{W}_h^i} \sigma_h^i(d_i u_h)(d_i w_h) = \sum_{i=1}^k \int_{\mathcal{K}_h^i} I_h(\sigma_h^i(d_i u_h)) d_i(\tilde{I}_h w_h).$$

Since $\mathcal{W}_h \subset \mathring{\mathcal{K}}_h$, $\tilde{I}_h w_h = 0$ on $\partial \mathcal{K}_h$, and then, performing an integration by part,

$$\sum_{i=1}^k \int_{\mathcal{W}_h^i} \sigma_h^i (d_i u_h)(d_i w_h) = -\sum_{i=1}^k \int_{\mathring{\mathcal{K}}_h} d_i \left( I_h(\sigma_h^i d_i u_h) \right) \tilde{I}_h(w_h)$$

$$= -\sum_{i=1}^k \int_{\mathcal{W}_h} d_i \left( I_h(\sigma_h^i d_i u_h) \right) w_h$$

$$= -\int_{\mathring{\mathcal{W}}_h} \sum_{i=1}^k d_i \left( I_h(\sigma_h^i d_i u_h) \right) w_h - \int_{\partial \mathcal{W}_h} \sum_{i=1}^k d_i \left( I_h(\sigma_h^i d_i u_h) \right) w_h.$$

To end the proof of Proposition 2.5, we remark that

$$I_h(\sigma_h^i d_i u_h) = \sigma_h^i d_i u_h \text{ on } \mathcal{W}_h^i$$

and then that

$$\sum_{i=1}^k d_i \left( I_h(\sigma^i d_i u_h) \right) = \Delta_h u_h \text{ on } \mathring{\mathcal{W}}_h.$$

Besides,

$$d_i \left( I_h(\sigma_h^i d_i u_h) \right) = \frac{-2}{h} a_i \left( I_h(\sigma_h^i d_i u_h) \right) n_h^i$$

and then

$$-\sum_{i=1}^k d_i \left( I_h(\sigma_h^i d_i u_h) \right) = \frac{1}{h} \partial_{n,h} u_h \text{ on } \partial \mathcal{W}_h,$$

thus concluding the proof of Proposition 2.5. □

To conclude this section, we also define the Sobolev norms for any trace function $g_h$ in $C(\partial \mathcal{W}_h)$ as:

$$|g_h|_{H^{1/2}(\partial \mathcal{W}_h)} = \min_{\substack{u_h \in \mathcal{W}_h \\ u_h = g_h \text{ on } \partial \mathcal{W}_h}} \|u_h\|_{H^1(\mathcal{W}_h)},$$

$$|g_h|_{H^{-1/2}(\partial \mathcal{W}_h)} = \max_{\substack{u_h \in \partial \mathcal{W}_h \\ |u_h|_{H^{1/2}(\partial \mathcal{W}_h)} = 1}} \frac{1}{h} \int_{\partial \mathcal{W}_h} g_h u_h.$$

Note that these definitions are of course compatible with the continuous ones.

### 2.4. The DtN map and the discrete Calderón problem

We are now in position to define the discrete Dirichlet-to-Neuman map (DtN map) for a potential $q_h \in C(\mathring{\mathcal{W}}_h)$:

**Definition 2.6** *Let $h > 0$. For any $q_h \in C(\mathring{\mathcal{W}}_h)$ such that*

$$\left( \Delta_h u_h + q_h u_h = 0 \text{ in } \mathring{\mathcal{W}}_h \text{ and } u_h = 0 \text{ on } \partial \mathcal{W}_h \right) \Rightarrow u_h = 0, \quad (2.17)$$

*we define the discrete Dirichlet-to-Neuman map (DtN map) as follows: $\Lambda_h[q_h]$ is an operator from $C(\partial \mathcal{W}_h)$ to $C(\partial \mathcal{W}_h)$ defined by $\Lambda_h[q_h](g_h) = \partial_{n,h} u_h$ where $\partial_{n,h} u_h$ is defined in (2.15) and $u_h$ is the unique solution in $C(\mathcal{W}_h)$ of (1.3) corresponding to the potential $q_h$.*

Note that the unique continuation property (2.17) is not very restrictive but is needed to guarantee the uniqueness of the solution $u_h$ of (1.3). First, the operator $\Delta_h + q_h$ is symmetric and derives from $(\Delta + q)$, a continuous operator with compact resolvent, hence its eigenspaces are finite dimensional (in the continuous setting), its eigenvalues are separated and a standard bootstrap argument shows that the elements in the kernel of $\Delta u + qu$ are regular (their second derivatives are as regular as $q$) hence if we suppose that the continuous operator $u \mapsto \Delta u + qu$ has a null kernel, any sensible numerical method will ensure that $u_h \mapsto \Delta_h u_h + q_h u_h$ has a null kernel for $h$ small enough.

Now, if the continuous potential $q$ derives from a conductivity $\sigma$ through Liouville transform through the formula $q = -\Delta(\sqrt{\sigma})/\sqrt{\sigma}$, automatically, solutions $u$ of (1.2) with $u = 0$ on the boundary correspond to solutions $v$ of

$$\text{div}(\sigma \nabla v) = 0 \text{ with } v = 0 \text{ on } \partial\Omega, \tag{2.18}$$

(more precisely, $u = \sqrt{\sigma} v$) and hence $v = u = 0$. Indeed, multiplying (2.18) by $v$ and integrating by parts, one immediately checks that $v \equiv 0$. This indicates that condition (2.17) is not very restrictive for practical applications.

The discrete Calderón problem then consists in the following one: Is the map $\Lambda_h$ injective?

But, as we have explained in the introduction, the main issue when dealing with a family of discrete Calderón problems is to get uniform stability estimates, where uniform means uniformly with respect to the mesh-size parameter $h > 0$.

To be more precise, we want to find a function $\omega : \mathbb{R}^+ \to \mathbb{R}^+$ such that, for all $h > 0$ and $(q_{1,h}, q_{2,h}) \in C(\mathring{\mathcal{W}}_h)^2$,

$$\|q_{1,h} - q_{2,h}\| \leq \omega\left(\|\Lambda_h[q_{1,h}] - \Lambda_h[q_{2,h}]\|\right).$$

Here, we do not give yet the norms that shall be put in the left and right hand-sides of this estimate, which actually contains much of the information.

Before stating our result precisely, let us introduce the discrete $H^r$-norms for $r \in \mathbb{R}$. For $h > 0$, set $\hat{\mathcal{K}}_h = [\![0, N-1]\!]^d$ and define the discrete Fourier transform of a function $u_h$ in $C_c(\mathcal{K}_h)$ as $\hat{u}_h$ in $C(\hat{\mathcal{K}}_h)$ by

$$\hat{u}_h(\xi) = \mathcal{F}_h(u_h)(\xi) = \int_{\mathcal{K}_h} u_h(x) e^{-2i\pi(x\cdot\xi)}, \quad \forall \xi \in \hat{\mathcal{K}}_h.$$

We then define the Sobolev norm $H^r$, for any $r \in \mathbb{R}$ by

$$|u_h|^2_{H^r(\mathcal{K}_h)} = \sum_{\xi \in \hat{\mathcal{K}}_h} |\hat{u}(\xi)|^2 (1 + |\xi|^2)^r.$$

Note that, with this definition, the celebrated Plancherel formula states that the norm $|u|_{H^0(\mathcal{K}_h)}$ coincides with $\|u\|_{L^2(\mathcal{K}_h)}$. Moreover, we emphasize, that, as in the continuous case, it is classical to show that $\|\cdot\|_{H^1(\mathcal{K}_h)}$ defined in (2.12) and $|\cdot|_{H^1(\mathcal{K}_h)}$ are equivalent independently of $h$.

Our main result is the following one:

**Theorem 2.7** *Suppose that Assumptions 1–4 are satisfied. Let $q_{1,h}$ and $q_{2,h}$ in $C(\mathring{\mathcal{W}}_h)$ be two potentials satisfying (2.17) and $\|q_{j,h}\|_{L^\infty(\mathring{\mathcal{W}}_h)} \leq m$, $j \in \{1,2\}$. Set*

$$\|\Lambda_h[q_{1,h}] - \Lambda_h[q_{2,h}]\|_{\mathfrak{L}_h} := \max_{\|g_h\|_{H^{1/2}(\partial \mathcal{W}_h)} = 1} |(\Lambda_h[q_{1,h}] - \Lambda_h[q_{2,h}])(g_h)|_{H^{-1/2}(\partial \mathcal{W}_h)}. \tag{2.19}$$

There exist constants $c, C$ depending on $m$ so that if $\|\Lambda_h[q_{1,h}] - \Lambda_h[q_{2,h}]\|_{\mathfrak{L}_h}, h, \epsilon_d$ and $\epsilon_a$ are smaller than $c$, setting

$$\mu := \max\left\{\epsilon_d, \epsilon_a^{1/2}, h^{1/2}, \frac{1}{|\log(\|\Lambda_h[q_{1,h}] - \Lambda_h[q_{2,h}]\|_{\mathfrak{L}_h})|}\right\}, \quad (2.20)$$

for all frequency $\xi \in \hat{\mathcal{K}}_h$ satisfying $|\xi| \leq c\mu^{-1}$,

$$|\mathcal{F}_h(q_{1,h} - q_{2,h})(\xi)| \leq C\mu. \quad (2.21)$$

Consequently, we have the following stability estimate : For any $r > 0$, there exists $C_r > 0$ depending on $r, m$ such that

$$|q_{1,h} - q_{2,h}|_{H^{-r}(\mathcal{W}_h)} \leq C_r \mu^{\frac{2r}{2r+d}}. \quad (2.22)$$

Before commenting this result, let us emphasize that the norm (2.19) corresponds to a discrete version of the $\mathfrak{L}(H^{1/2}(\partial\Omega); H^{-1/2}(\partial\Omega))$-norm, thus being completely consistent with the natural norms for the continuous Calderón problem.

Also note that Theorem 2.7 actually bears four results.

- The influence of the mesh size alone is seen by setting $\mu = h^{1/2}$ and $\Lambda_h[q_{1,h}] = \Lambda_h[q_{2,h}]$, this is the case, for instance for a regular mesh where the DtN operators are equal, the uncertainty on the potentials $q_h$ in the $H^{-r}(\mathcal{W}_h)$ norm then scales as $h^{\frac{r}{2r+d}}$.

- The influence of the anisotropic scaling $\epsilon_a$ or the irregularity of the mesh $\epsilon_d$ is measured by taking $\mu = \epsilon_d$ or $\epsilon_a^{1/2}$.

- The influence of the error in the DtN maps, this is the true stability estimate. The stability estimate is in $|\log(error)|^\alpha$, with $\alpha < 0$. It is consistent with every continuous stability estimate (see [1]).

At this point, let us also emphasize that Theorem 2.7 does not yield uniqueness. Indeed, if $\Lambda_h[q_{1,h}] = \Lambda_h[q_{2,h}]$, then for any $r > 0$, we only have

$$|q_{1,h} - q_{2,h}|_{H^{-r}(\mathcal{W}_h)} \leq C_r \max\{\epsilon_d, \epsilon_a^{1/2}, h^{1/2}\}^{\frac{2r}{2r+d}}. \quad (2.23)$$

Of course, given $h > 0$, this does not give uniqueness, but rather some kind of asymptotic uniqueness as $h \to 0$, since the right hand-side of (2.23) goes to zero as $h \to 0$.

The proof of Theorem 2.7 is the main goal of that article. It will be given in Section 5 and it will strongly use the results developed in the other sections, and in particular the construction of discrete CGO solutions derived in Section 4.

### 2.5. Basic properties of the DtN maps

Let us begin our analysis of the discrete Calderón problems by giving some properties of the discrete DtN problems that are classical in the continuous setting:

**Proposition 2.8** *Let $h > 0$ and $q_h \in C(\mathring{\mathcal{W}}_h)$ satisfying (2.17).*

*The discrete DtN map is self adjoint:*

$$\int_{\partial \mathcal{W}_h} \Lambda_h[q_h](u_h) v_h = \int_{\partial \mathcal{W}_h} \Lambda_h[q_h](v_h) u_h, \quad \forall (u_h, v_h) \in C(\partial \mathcal{W}_h)^2.$$

*Moreover, for all $u_{1,h}$ and $u_{2,h}$ that verify $\Delta_h u_{j,h} + q_{j,h} u_{j,h} = 0$ for $j = 1, 2$ on $\mathring{\mathcal{W}}_h$ with $q_{1,h}, q_{2,h} \in C(\mathring{\mathcal{W}}_h)$ satisfying (2.17), we have*

$$\int_{\mathring{\mathcal{W}}_h} (q_{1,h} - q_{2,h}) u_{1,h} u_{2,h} = \frac{1}{h} \int_{\partial \mathcal{W}_h} (\Lambda_h[q_{1,h}] - \Lambda_h[q_{2,h}])(u_{1,h}) u_{2,h}, \quad (2.24)$$

where we used the slight abuse of notations consisting of denoting the same way a function on $\mathcal{W}_h$ and its trace on the boundary $\partial \mathcal{W}_h$.

The proof of these assertions is similar to the one of the continuous case.

**Proof** We first prove the self-adjointness of $\Lambda_h[q_h]$ for $q_h \in C(\mathring{\mathcal{W}}_h)$. For $(u_h, v_h) \in C(\mathcal{W}_h)^2$, solutions of $\Delta_h u_h + q_h u_h = 0 = \Delta_h v_h + q_h v_h$ on $\mathring{\mathcal{W}}_h$, using the integration by part equation (2.16), we obtain

$$\frac{1}{h} \int_{\partial \mathcal{W}_h} \Lambda_h[q_h](u_h) v_h = \frac{1}{h} \int_{\partial \mathcal{W}_h} (\partial_{n,h} u_h) v_h$$
$$= \int_{\mathring{\mathcal{W}}_h} (\Delta_h u_h) v_h + \sum_{i=1}^{k} \int_{\mathcal{W}_h^i} \sigma_h^i (d_i u_h)(d_i v_h)$$
$$= \int_{\mathring{\mathcal{W}}_h} -q_h u_h v_h + \sum_{i=1}^{k} \int_{\mathcal{W}_h^i} \sigma_h^i (d_i u_h)(d_i v_h).$$

The right-hand side being symmetric in $u_h$ and $v_h$, their roles can be exchanged which proves the self adjointness of $\Lambda_h[q_h]$.

We now turn our attention to (2.24). Set $u_{1,h}$ and $u_{2,h}$ as in Proposition 2.8, then

$$\int_{\mathring{\mathcal{W}}_h} q_{1,h} u_{1,h} u_{2,h} = -\int_{\mathring{\mathcal{W}}_h} (\Delta_h u_{1,h}) u_{2,h}$$
$$= -\frac{1}{h} \int_{\partial \mathcal{W}_h} (\partial_{n,h} u_{1,h}) u_{2,h} + \sum_{i=1}^{k} \int_{\mathcal{W}_h^i} \sigma_h^i (d_i u_{1,h})(d_i u_{2,h}).$$

Recalling that, by definition $\partial_{n,h} u_{1,h} = \Lambda_h[q_{1,h}](u_{1,h})$, interchanging the roles of $u_{1,h}$ and $u_{2,h}$ and subtracting the two equations lead to

$$\int_{\mathring{\mathcal{W}}_h} (q_{1,h} - q_{2,h}) u_{1,h} u_{2,h} = -\frac{1}{h} \int_{\partial \mathcal{W}_h} (\Lambda_h[q_{1,h}](u_{1,h}) u_{2,h} - \Lambda_h[q_{2,h}](u_{2,h}) u_{1,h})$$
$$= -\frac{1}{h} \int_{\partial \mathcal{W}_h} (\Lambda_h[q_{1,h}] - \Lambda_h[q_{2,h}])(u_{1,h}) u_{2,h},$$

where we used the self-adjointness of $\Lambda_h[q_{2,h}]$ in order to conclude. □

Note that formula (2.24) will be the basis for the proof of the stability estimates in Theorem 2.7, as it will be seen in Section 5.

### 2.6. Notations

In the sequel, we shall denote by $C$ and $c$ generic constants that can be chosen independently of the parameters $h$, $\epsilon_a$, $\epsilon_d$, the Carleman parameter $s$ and the CGO parameter $\beta, s, \eta$. The constant $C$ will be chosen large enough, whereas $c$ will be a small positive constant.

## 3. Carleman estimate

The main tool that proves the existence of CGO solutions is a Carleman estimate proved with special weight functions called the "limiting Carleman weights" (see [18]). We emphasize the fact that discrete Carleman estimates already exist, see [6, 7], but that they are not based on limiting Carleman weights.

The goal of this section is to develop the proof of uniform discrete Carleman estimate (here, uniform means with respect to the discretization parameter $h > 0$) with a simple limiting Carleman weight, namely the one corresponding to plane waves.

*3.1. Statement of the Carleman estimate*

For any $s \in \mathbb{R}^d$, introduce $\phi_s(x) = s \cdot x$ and define the operator $\Delta_{s,h}$ from $C_c(\mathcal{K}_h)$ to $C_c(\mathcal{K}_h)$ as

$$\Delta_{s,h} u_h = \begin{cases} \sum_{i=1}^k e^{-\phi_s} d_i \left( \sigma_h^i d_i (e^{\phi_s} u_h) \right) & \text{on } \mathring{\mathcal{K}}_h \\ 0 \text{ on } \partial \mathcal{K}_h \end{cases}, \quad (3.1)$$

We aim at performing a Carleman estimate on $\Delta_{s,h}$.

Below, we use the slight abuse of notation consisting of identifying functions of $C_c(\mathcal{B}_h)$ as functions of $C_c(\mathcal{K}_h)$, by extending functions of $C_c(\mathcal{B}_h)$ by zero outside $\mathcal{B}_h$.

**Theorem 3.1** *There exist $c > 0$ and $C > 0$ such that for all $(h, \epsilon_d, \epsilon_a) \in (0, c)^3$, for all $s \in \mathbb{R}^d$ such that $|s| \leq c \min\{\epsilon_d^{-1}, h^{-2/3}\}$, for any $u_h \in C_c(\mathcal{B}_h)$, we have*

$$|s| \|u_h\|_{L^2(\mathcal{B}_h)} + \|u_h\|_{\dot{H}^1(\mathcal{B}_h)} \leq C \|\Delta_{s,h} u_h\|_{L^2(\mathcal{K}_h)}. \quad (3.2)$$

Note that estimate (3.2) is similar to the continuous one. However, an important difference between (3.2) and the continuous Carleman estimate is that the range of $s$ for which (3.2) holds true is limited to some scales depending on the mesh. Of course, these scales are going to infinity as $(h, \epsilon_d) \to 0$, thus being completely compatible with the continuous Carleman estimate.

**Remark 3.2** *When $k = d$ (i.e, there are as many connections as the dimension), then the scaling becomes $|s| \leq c \min\{\epsilon_d^{-1}, h^{-1}\}$ and thus we gain an order of $s$ with respect to $h$. This scaling is coherent with the one obtained in [7]. We do not state it as a main result since the final scaling appearing in the stability theorem (Theorem 2.7) is $|s| \leq h^{-1/2}$ and is driven by Proposition 4.1. Nevertheless, we show in Subsection 3.6 how to improve the following proof in order to obtain the scaling $h^{-1}$.*

**Proof** [Sketch of the proof] The proof of (3.2) is based on three steps that will be developed in details in next sections:

- Decompose the operator $\Delta_{s,h}$ into its symmetric and skew-symmetric parts, $S_{s,h}$ and $A_{s,h}$ (see Proposition 3.3);
- Prove a lower bound on $\|S_{s,h} u_h\|_{L^2(\mathcal{K}_h)} + \|A_{s,h} u_h\|_{L^2(\mathcal{K}_h)}$ of the form

$$|s| \|u_h\|_{L^2(\mathcal{B}_h)} + \|u_h\|_{\dot{H}^1(\mathcal{B}_h)} \leq C \left( \|S_{s,h} u_h\|_{L^2(\mathcal{K}_h)} + \|A_{s,h} u_h\|_{L^2(\mathcal{K}_h)} \right),$$

see Proposition 3.4 for more precise statements;
- Prove an upper bound on $\int_{\mathcal{K}_h} (S_{s,h} u_h)(A_{s,h} u_h)$ of the form

$$\int_{\mathcal{K}_h} (S_{s,h} u_h)(A_{s,h} u_h) = o(\|S_{s,h} u_h\|_{L^2(\mathcal{K}_h)}^2 + \|A_{s,h} u_h\|_{L^2(\mathcal{K}_h)}^2)$$

for all $u_h \in C_c(\mathcal{B}_h)$, see Proposition 3.9 for detailed statements.

Then the decomposition

$$\|\Delta_{s,h} u_h\|_{L^2}^2 = \|S_{s,h} u_h\|_{L^2}^2 + \|A_{s,h} u_h\|_{L^2}^2 + 2 \int_{\mathcal{K}_h} (S_{s,h} u_h)(A_{s,h} u_h)$$

yields (3.2) immediately, see Subsection 3.5 for the final compilation of the different listed propositions. □

Also note that actually, for Theorem 3.1 to be true, we strongly need $\mathcal{B}_h$ to be such that $\mathcal{B}_h \subset \mathring{\mathcal{K}}_h$, thus explaining where this condition comes from in Assumption 4. The reason is we need to perform several integration by parts without boundary terms in the proof of the upper bound of the commutator. We then need to ensure that the extension by 0 of $u_h$ outside $\mathcal{B}_h$ satisfies that $u_h$, $d_i u_h$ and $a_j d_i u_h$ equal to 0 on, respectively, $\partial \mathcal{K}_h$, $\partial(\mathcal{K}_h^i)$, $\partial(\mathcal{K}_h^{ij})$.

To make easier the reading of the tedious computations that will come afterwards, let us briefly recall how the Carleman estimate is proved in the continuous case. First, one easily checks that

$$\Delta_s u = e^{-\phi_s} \Delta(e^{\phi_s} u) = \Delta u + 2s \cdot \nabla u + |s|^2 u.$$

This operator is then decomposed into its symmetric and skew-symmetric parts:

$$S_s u = \Delta u + |s|^2 u \text{ and } A_s u = 2s \cdot \nabla u. \tag{3.3}$$

Then Poincaré inequality yields

$$\|A_s u\|_{L^2(\Omega)} \gtrsim |s| \|u\|_{L^2(\Omega)}, \tag{3.4}$$

whereas

$$\int_\Omega |\nabla u|^2 = -\int_\Omega (S_s u) u + |s|^2 \int_\Omega |u|^2 \leq \|S_s u\|_{L^2} \|u\|_{L^2} + |s|^2 \|u\|_{L^2}.$$

Using the previous estimate, we thus directly get

$$\|u\|_{\dot{H}^1(\Omega)}^2 + |s|^2 \|u\|_{L^2(\Omega)}^2 \lesssim \|S_s u\|_{L^2(\Omega)}^2 + \|A_s u\|_{L^2(\Omega)}^2,$$

which concludes the proof in the continuous case since one easily checks that $\int_\Omega S_s u A_s u = 0$ for $u \in C_c^\infty(\Omega)$.

*3.2. Decomposition into the symmetric and skew-symmetric part*

The discrete counterpart of the continuous decomposition in (3.3) is:

**Proposition 3.3** *Let $h > 0$. For $i \in \{1, \cdots, k\}$, define $\kappa_{s,h}^i, g_{s,h}^i, p_{s,h}^i \in \mathbb{R}$ as*

$$\kappa_{s,h}^i := \text{ch}(hs \cdot e_i), \quad p_{s,h}^i := \frac{4}{h^2} \text{sh}^2\left(\frac{h}{2} s \cdot e_i\right), \quad g_{s,h}^i := \frac{1}{h} \text{sh}(hs \cdot e_i). \tag{3.5}$$

*Define the operators $S_{s,h}^i$ and $A_{s,h}^i$ from $C_c(\mathcal{K}_h)$ to $C_c(\mathcal{K}_h)$ as*

$$S_{s,h}^i u_h := \kappa_{s,h}^i d_i(\sigma_h^i d_i u_h) + p_{s,h}^i (a_i \sigma_h^i) u_h \text{ on } \mathring{\mathcal{K}}_h$$
$$A_{s,h}^i u_h := 2 g_{s,h}^i a_i(\sigma_h^i d_i u_h) + g_{s,h}^i (d_i \sigma_h^i) u_h \text{ on } \mathring{\mathcal{K}}_h.$$

*Then the $S_{s,h}^i$ are symmetric whereas the $A_{s,h}^i$ are skew-symmetric and we have*

$$\Delta_{s,h} = S_{s,h} + A_{s,h} \quad \text{with } S_{s,h} = \sum_{i=1}^k S_{s,h}^i \text{ and } A = \sum_{i=1}^k A_{s,h}^i.$$

**Proof** Before going into the proof, let us emphasize that it can be done for $s \in \mathbb{R}^d$ and $h > 0$ fixed. We shall therefore omit the indexes $s$ and $h$ in the proof to simplify the notations.

Note that thanks to (2.3) in Lemma 2.1, we have

$$a_i(e^\phi) d_i(e^{-\phi}) = -a_i(e^{-\phi}) d_i(e^\phi)$$

The proof of Proposition 3.3 is then an application of Lemma 2.1 and a discrete integration by parts (see Proposition 2.4) with $v = 0$ on $\partial \mathcal{K}_h$. For all $u, v \in C_c(\mathcal{K}_h)$, we have:

$$-\int_{\mathcal{K}_h} (\Delta_{s,h} u) v = -\int_{\mathring{\mathcal{K}}_h} \left( \sum_i e^{-\phi} d_i \left( \sigma^i d_i (e^\phi u) \right) \right) v$$

$$= \sum_i \int_{\mathcal{K}_h^i} \sigma^i d_i(u e^\phi) d_i(v e^{-\phi})$$

$$= \sum_i \int_{\mathcal{K}_h^i} \sigma^i \left( d_i(u) a_i(e^\phi) + d_i(e^\phi) a_i(u) \right) \left( d_i(v) a_i(e^{-\phi}) + d_i(e^{-\phi}) a_i(v) \right)$$

$$= \sum_i \int_{\mathcal{K}_h^i} \left( r_1^i d_i(u) d_i(v) + r_2^i a_i(u) a_i(v) + r_3^i d_i(u) a_i(v) - r_3^i a_i(u) d_i(v) \right)$$

where a direct computation, using Lemma 2.2, yields :

$$r_1^i = \sigma^i a_i(e^\phi) a_i(e^{-\phi}) = \sigma^i ch^2 \left( \frac{h}{2} s \cdot e_i \right), \tag{3.6}$$

$$r_2^i = \sigma^i d_i(e^\phi) d_i(e^{-\phi}) = \frac{-4\sigma^i}{h^2} sh^2 \left( \frac{h}{2} s \cdot e_i \right), \tag{3.7}$$

$$r_3^i = \sigma^i a_i(e^\phi) d_i(e^{-\phi}) = -\frac{\sigma^i}{h} sh(hs \cdot e_i). \tag{3.8}$$

Hence, for all $i \in \{1, \cdots, k\}$ the decomposition into the symmetric part and the skew-symmetric part of the operator $e^{-\phi} d_i(\sigma^i d_i(e^\phi u))$ is given by

$$(S_i u, v) = -\int_{\mathcal{K}_h^i} r_1^i d_i(u) d_i(v) + r_2^i a_i(u) a_i(v)$$

$$= \int_{\mathring{\mathcal{K}}_h} [d_i(r_1^i d_i u) - a_i(r_2^i a_i u)] v$$

$$(A_i u, v) = -\int_{\mathcal{K}_h^i} r_3^i (d_i(u) a_i(v) - d_i(v) a_i(u))$$

$$= -\int_{\mathring{\mathcal{K}}_h} [a_i(r_3^i d_i u) + d_i(r_3^i a_i u)] v.$$

Using Lemma 2.3, we obtain, on $\mathring{\mathcal{W}}_h$ :

$$S_i = d_i \left( r_1^i - \frac{h^2}{4} r_2^i \right) d_i u - (a_i r_2^i) u, \quad A_i = -2 a_i (r_3^i d_i u) - (d_i r_3^i) u.$$

Then, basic identities on the hyperbolic functions yield Proposition 3.3. □

*3.3. The Lower bound*

The goal of this section is to prove the following proposition:

**Proposition 3.4** *Let $S_{s,h}, A_{s,h}$ be as in Proposition 3.3. Then there exist $c > 0$ and $C > 0$ such that for $(h, \epsilon_d, \epsilon_a) \in (0, c)^3$, for all $|s| \leq ch^{-2/3}$, for all $u_h \in C_c(\mathcal{B}_h)$,*

$$|s| \|u_h\|_{L^2(\mathcal{K}_h)} + \|u_h\|_{\dot{H}^1(\mathcal{K}_h)} \leq C (\|S_{s,h} u_h\|_{L^2(\mathcal{K}_h)} + \|A_{s,h} u_h\|_{L^2(\mathcal{K}_h)}) \tag{3.9}$$

The proof of Proposition 3.4 is postponed to the end of the section.

As in the continuous case, see above, the first step of such a proof is a Poincaré inequality that yields (3.4).

**Lemma 3.5** *Given $s \in \mathbb{R}^d$ and $h > 0$, define the operator $\hat{A}_{s,h}$ on $C_c(\mathcal{K}_h)$ by*

$$\hat{A}_{s,h} u_h := 2 \sum_{i=1}^{k} (s \cdot e_i)(a_i \sigma_h^i) a_i d_i u_h. \tag{3.10}$$

*Then there exist $c > 0$ and $C > 0$ such that for all $(\epsilon_d, \epsilon_a) \in (0, c)^2$, for all $s \in \mathbb{R}^d$,*

$$|s|^2 \|u_h\|_{L^2(\mathcal{K}_h)}^2 \leq C \|\hat{A}_{s,h} u_h\|_{L^2(\mathcal{K}_h)}^2, \quad \forall u_h \in C_c(\mathcal{B}_h). \tag{3.11}$$

**Proof** Below, we use the symbol $\simeq$ to refer to an equality modulo a term of order $(\epsilon_d + \epsilon_a)|s|^2 \|u\|_{L^2(\mathcal{K}_h)}^2$. By setting $c$ small enough these terms are negligible when compared to $|s|^2 \|u_h\|_{L^2(\mathcal{K}_h)}^2$.

$$|s|^2 \|u\|_{L^2(\mathcal{K}_h)}^2 \simeq \int_{\mathcal{K}_h} \left( \sum_{i=1}^{k} (s \cdot e_i)^2 a_i(\sigma^i) \right) u^2$$

$$= \int_{\mathcal{K}_h} \sum_{i=1}^{k} a_i d_i(s \cdot x) a_i(\sigma^i)(s \cdot e_i) u^2$$

$$= -\int_{\mathcal{K}_h} \sum_{i=1}^{k} (s \cdot x) a_i d_i \left( (s \cdot e_i) a_i(\sigma^i) u^2 \right)$$

$$= -\int_{\mathcal{K}_h} (s \cdot x) \sum_{i=1}^{k} (s \cdot e_i) a_i d_i \left( a_i(\sigma^i) u^2 \right)$$

$$\simeq -\int_{\mathcal{K}_h} (s \cdot x) \sum_{i=1}^{k} (s \cdot e_i) a_i(\sigma^i) a_i d_i \left( u^2 \right)$$

$$= -2 \int_{\mathcal{K}_h} (s \cdot x) \sum_{i=1}^{k} (s \cdot e_i) a_i(\sigma^i) a_i d_i(u) \left( a_i a_i(u) + \frac{h^2}{4} d_i d_i(u) \right)$$

$$\leq C |s| \|\hat{A} u\|_{L^2(\mathcal{K}_h)} \|u\|_{L^2(\mathcal{K}_h)}.$$

Here, for the first sign $\simeq$, we used (2.8), for the second one we used $|d_i(\sigma^i)| \leq \epsilon_d$ (see (2.7)) and, for the last estimate, that $\|a_i a_i(u)\| + h^2 \|d_i d_i u\| \leq C \|u\|$. $\square$

**Remark 3.6** *Note that this discrete Poincaré inequality uses dearly the fact $u_h = 0$ on $\partial \mathcal{B}_h$ and that $\mathcal{B}_h \subset \ddot{\mathcal{K}}_h$. This is not only a technical argument for the proof. Indeed $u_h$ must cancel on the boundary of $\mathcal{K}_h$ (as expected for any Poincaré inequality) but also on the points at distance one of the boundary.*

*We can indeed build a counterexample to the Poincare inequality in (3.11) in 1-d, by taking $u_h(jh) = 0$ for $j$ even and $u_h(jh) = 1$ for $j$ odd: in this case, $ad(u_h) = 0$ whereas $\|u_h\| \neq 0$ and $u_h$ indeed cancels on the boundary if $\mathcal{B}_h = \{0, \cdots, jh, \cdots, Nh\}$ for $N$ even.*

We are now in position to prove Proposition 3.4.

**Proof** [Proof of Proposition 3.4] We now turn our attention to proving (3.9). First we remark that $\kappa^i \geq 1$ (defined in (3.5)) and $|p^i a_i(\sigma^i)| \leq C|s|^2$ (since $|s|h$ is uniformly bounded), hence we have, for all $u \in C_c(\mathcal{B}_h)$,

$$-(Su, u) \geq \frac{1}{2} \int_{\mathcal{K}_h} \sum_i \sigma^i d_i(u) d_i(u) - C|s|^2 \|u\|_{L^2(\mathcal{K}_h)}^2$$

$$= \frac{1}{2} \|u\|_{\dot{H}^1(\mathcal{K}_h)}^2 - C|s|^2 \|u\|_{L^2(\mathcal{K}_h)}^2.$$

Using
$$|(Su,u)| \leq \|Su\|\|u\| \leq \frac{1}{2}(\|Su\|^2 + \|u\|^2),$$
we get
$$\|u\|^2_{\dot{H}^1(\mathcal{K}_h)} \leq C\left(\|Su\|^2_{L^2(\mathcal{K}_h)} + |s|^2\|u\|^2_{L^2(\mathcal{K}_h)}\right).$$

Now, Lemma 3.5 yields
$$|s|^2\|u\|^2_{L^2(\mathcal{K}_h)} + \|u\|^2_{\dot{H}^1(\mathcal{K}_h)} \leq C\left(\|Su\|^2_{L^2(\mathcal{K}_h)} + \|\hat{A}_{s,h}u\|^2_{L^2(\mathcal{K}_h)}\right), \quad (3.12)$$
where $\hat{A}_{s,h}$ is the operator defined in (3.10).

To conclude, we shall then compare the operators $\hat{A}_{s,h}$ and $A_{s,h}$.

Since $|g^i - s \cdot e_i| \leq C|s|^3 h^2$ (recall the definition of $g^i$ in (3.5) and $|s|h$ uniformly bounded), for all $u \in C_c(\mathcal{B}_h)$,
$$\|2\sum_i g^i(a_i\sigma^i) a_i d_i u - \hat{A}u\|_{L^2(\mathcal{K}_h)} \leq C|s|^3 h^2 \|u\|_{\dot{H}^1(\mathcal{K}_h)},$$
and thus, for every $|s| \leq ch^{-2/3}$,
$$\|2\sum_i g^i(a_i\sigma^i) a_i d_i u - \hat{A}u\|^2_{L^2(\mathcal{K}_h)} \leq cC\|u\|^2_{\dot{H}^1(\mathcal{K}_h)}. \quad (3.13)$$

Besides, using that $|g^i| \leq C|s|$ and $|d_i\sigma^i| \leq \epsilon_d$ by (2.7),
$$\|g^i(d_i\sigma^i)u\|^2_{L^2(\mathcal{K}_h)} \leq C\epsilon_d^2 |s|^2 \|u\|^2_{L^2(\mathcal{K}_h)}. \quad (3.14)$$

Using the two estimates (3.13)–(3.14) and plugging them in (3.12) while using the definition of $A^i$ in Proposition 3.3, we obtain (3.9) for $c$ small enough. □

**Remark 3.7** *The case* $k = d$. Let us now turn our attention to the case $k = d$. A careful reading of the previous proof shows that the only occurrence of the hypothesis $|s| \leq ch^{-2/3}$ lies in the proof of (3.13). In order to avoid this hypothesis, when $k = d$, we can define the vector $G$ such that $G \cdot e_i = g^i_{s,h}$, since $e_i$ form a basis of $\mathbb{R}^d$. Thus following the proof of Lemma 3.5, we can replace the operator $\hat{A}$ by $\tilde{A}$ defined as
$$\tilde{A} := 2\sum_{i=1}^k (G \cdot e_i)(a_i\sigma^i_h) a_i d_i u_h$$
and obtain, for $\epsilon_d + \epsilon_a$ small enough,
$$|G|^2 \|u_h\|^2_{L^2(\mathcal{K}_h)} \leq C\|\tilde{A}_{s,h} u_h\|^2_{L^2(\mathcal{K}_h)}, \quad \forall u_h \in C_c(\mathcal{B}_h). \quad (3.15)$$

Hence, since $\tilde{A}$ satisfies, similarly as in (3.14),
$$\|\tilde{A}u - Au\|^2_{L^2(\mathcal{K}_h)} \leq C\epsilon_d^2 |s|^2 \|u\|^2_{L^2(\mathcal{K}_h)},$$
following the proof of Proposition 3.4, one can prove:

**Proposition 3.8** *When* $k = d$, *there exist* $c > 0$ *and* $C > 0$ *such that for* $(h, \epsilon_d, \epsilon_a) \in (0, c)^3$, *for all* $|s| \leq ch^{-1}$, *for all* $u_h \in C_c(\mathcal{B}_h)$,
$$|s|\|u_h\|_{L^2(\mathcal{K}_h)} + \|u_h\|_{\dot{H}^1(\mathcal{K}_h)} \leq C(\|S_{s,h} u_h\|_{L^2(\mathcal{K}_h)} + \|A_{s,h} u_h\|_{L^2(\mathcal{K}_h)}). \quad (3.16)$$

Where we used $|s| \leq |g| \leq C_0|G|$, where $C_0$ depends only on the basis $(e_i)_i$.

*3.4. The upper bound of the commutator*

**Proposition 3.9** *Let $S_{s,h}, A_{s,h}$ be as in Proposition 3.3. There exists $C, c > 0$ such that for all $(h, \epsilon_a, \epsilon_d) \in (0, c)^3$, for all $\varepsilon > 0$, there exists $c_\varepsilon > 0$, depending on $\varepsilon$ such that for all $s \in \mathbb{R}^d$ such that $|s| \leq c_\varepsilon \min\{\epsilon_d^{-1}, h^{-2/3}\}$, then for all $u_h \in C_c(\mathcal{B}_h)$,*

$$\left| \int_{\mathring{\mathcal{K}}_h} (A_{s,h} u_h)(S_{s,h} u_h) \right| \leq \varepsilon (\|A_{s,h} u_h\|_{L^2(\mathcal{K}_h)}^2 + \|S_{s,h} u_h\|_{L^2(\mathcal{K}_h)}^2). \quad (3.17)$$

**Remark 3.10** *In the case $k = d$, following the proof of Proposition 3.9 and using Proposition 3.8, one can prove (3.17) for $|s| \leq c_\varepsilon \min\{\epsilon_d^{-1}, h^{-1}\}$.*

**Proof** First suppose that $c_\varepsilon < c$ so that Proposition 3.4 holds. As before, for simplicity of notations, in the proof below, we shall omit the indexes $s$ and $h$ since the proof is done for $s$ and $h$ fixed.

We shall introduce the slightly modified operators $\tilde{S}_i$ and $\tilde{A}_i$ defined on $C_c(\mathcal{K}_h)$ by

$$\tilde{S}_i u := (\kappa^i a_i(\sigma^i)) d_i d_i u + (p^i a_i(\sigma^i)) u = S_i u - (\kappa^i d_i(\sigma^i)) a_i d_i u$$

$$\tilde{A}_i u := 2(g^i a_i(\sigma^i)) a_i d_i u = A_i u - \frac{h^2}{4} g^i d_i(\sigma^i)) d_i d_i u - (g^i d_i(\sigma^i)) u$$

Using the bounds $\|u\|_{\dot{H}^1(\mathcal{K}_h)} \leq C(\|Su\|_{L^2(\mathcal{K}_h)} + \|Au\|_{L^2(\mathcal{K}_h)})$, $\kappa^i \leq C$ (since $|s|h$ is bounded) and $|d_i(\sigma)| \leq \epsilon_d$, then

$$\|\tilde{S}_i u - S_i u\|_{L^2(\mathcal{K}_h)} \leq \epsilon_d (\|Au\|_{L^2(\mathcal{K}_h)} + \|Su\|_{L^2(\mathcal{K}_h)}).$$

A similar development for the operators $A_i$ may be obtained by using

$$h^2 \|d_i d_i u\|_{L^2(\mathcal{K}_h)} \leq C \|u\|_{L^2(\mathcal{K}_h)} \leq C |s|^{-1} (\|Su\|_{L^2(\mathcal{K}_h)} + \|Au\|_{L^2(\mathcal{K}_h)}),$$

by using $|g^i| \leq C |s|$ (since $|s|h$ is bounded) and $|d(\sigma^i)| \leq \epsilon_d$ by (2.7).

Therefore, it is enough to prove Proposition 3.9 while replacing $S$ and $A$ by $\tilde{S} = \sum_i \tilde{S}_i$ and $\tilde{A} = \sum_i \tilde{A}_i$.

Define $L^{ij} \in C(\mathring{\mathcal{K}}_h)$

$$L^{ij} = L^{ji} := a_i(\sigma^i) a_j(\sigma^j) \quad \text{on } \mathring{\mathcal{K}}_h.$$

We then have, for all $u \in C_c(\mathcal{K}_h)$,

$$\frac{1}{2}(\tilde{S}_i u, \tilde{A}_j u) = (\kappa^i g^j) \int_{\mathring{\mathcal{K}}_h} L^{ij} (d_i d_i u)(a_j d_j u) + (p^i g^j) \int_{\mathring{\mathcal{K}}_h} L^{ij} u (a_j d_j u).$$

From the explicit form of the coefficients in (3.5), $\kappa^i \leq C$, $g^i \leq C|s|$ and $p^i \leq C|s|^2$. Therefore, since we shall ensure that $\epsilon_d |s|$ can be made arbitrary small at the end, it is sufficient to show the following lemma:

**Lemma 3.11** *For all $(i, j) \in \{1, \cdots, k\}^2$ and $s \in \mathbb{R}^d$, there exists a constant $C$ independent of $s$ and $\epsilon_d$ such that for all $u_h \in C_c(\mathcal{B}_h)$,*

$$\left| \int_{\mathring{\mathcal{K}}_h} L^{ij} (d_i d_i u)(a_j d_j u) \right| \leq C \epsilon_d (\|Su\|_{L^2(\mathcal{K}_h)}^2 + \|Au\|_{L^2(\mathcal{K}_h)}^2), \quad (3.18)$$

*and*

$$\left| \int_{\mathring{\mathcal{K}}_h} L^{ij} u (a_j d_j u) \right| \leq C |s|^{-2} \epsilon_d (\|Su\|_{L^2(\mathcal{K}_h)}^2 + \|Au\|_{L^2(\mathcal{K}_h)}^2). \quad (3.19)$$

These two estimates are established in Section 3.4.1 and 3.4.2 respectively and are based on two main ingredients:

- Thanks to the bound $\|d_i(\sigma^i)\|_{L^\infty} \le \epsilon_d$, we have $\|d_i(L^{ij})\|_{L^\infty} \le CM\epsilon_d$.
- $\mathcal{B}_h \subset \ddot{\mathcal{K}}_h$ so that for all $u \in C_c(\mathcal{B}_h)$

$$a_j u = d_j u = 0 \text{ on } \partial(\mathcal{K}_h^j) \text{ and } d_j d_j u = a_j d_j u = a_j a_j u = 0 \text{ on } \partial\mathcal{K}_h,$$

thus allowing us to perform several integration by parts without any boundary term.

The proofs of these estimates are given in the next sections. □

*3.4.1. Proof of estimate (3.18) in Lemma 3.11* Our proof is based on several integration by parts and interchanges between the operators $a$ and $d$.

The integrations by part are done in the following order: first on $a_i$, then on $d_i$, $d_j$ and finally on the last $d_i$. The underlined terms will scale as $C\epsilon_d(\|Su\|^2_{L^2(\mathcal{K}_h)} + \|Au\|^2_{L^2(\mathcal{K}_h)})$ and thus will be removed. They will be of the form:

$$d(L)d(u)d(u) \text{ (first case) or } h^2 d(L)d^3(u)d(u) \text{ (second case)}.$$

Terms corresponding to the first case are underlined before being removed. Note that those terms have the required scaling since $\|d(L^{ij})\|_{L^\infty} \le C\epsilon_d$ and $\|u\|_{\dot{H}^1(\mathcal{B}_h)} \le C(\|Su\|_{L^2(\mathcal{K}_h)} + \|Au\|_{L^2(\mathcal{K}_h)})$ and $\|h^2 d^2(v)\| \le C\|v\|$. Below, the symbol $\simeq$ will refer to an equality modulo a term $C\epsilon_d(\|Su\|^2_{\mathcal{K}_h} + \|Au\|^2_{\mathcal{K}_h})$.

$$\int_{\mathring{\mathcal{K}}_h} (d_i d_i u) L^{ij} (a_j d_j u) = \int_{\mathring{\mathcal{K}}_h} (d_i d_i u) L^{ij} (a_j d_j u)$$

$$= \int_{\mathcal{K}_h^j} a_j [(d_i d_i u) L^{ij}](d_j u) \quad (\text{since } d_i d_i u = 0 \text{ on } \partial \mathcal{K}_h)$$

$$= \int_{\mathcal{K}_h^j} (a_j d_i d_i u)(a_j L^{ij})(d_j u) + \underbrace{\frac{h^2}{4} \int_{\mathcal{B}_h^j} (d_j d_i d_i u)(d_j L^{ij})(d_j u)}_{\text{2nd case}}$$

$$\simeq - \int_{\mathcal{K}_h^{ij}} (a_j d_i u) d_i[(a_j L^{ij})(d_j u)] \quad (\text{since } d_j u = 0 \text{ on } \partial(\mathcal{K}_h^j))$$

$$= - \int_{\mathcal{K}_h^{ij}} (a_j d_i u)(a_i a_j L^{ij})(d_i d_j u) - \underbrace{\int_{\mathcal{K}_h^{ij}} (a_j d_i u)(d_i a_j L^{ij})(a_i d_j u)}_{\text{1st case}}$$

$$\simeq \int_{\mathcal{K}_h^i} d_j[(a_j d_i u)(a_i a_j L^{ij})](d_i u) \quad (\text{since } d_i u = 0 \text{ on } \partial(\mathcal{K}_h^i))$$

$$= \int_{\mathcal{K}_h^i} (d_j a_j d_i u)(a_j a_i a_j L^{ij})(d_i u) + \underbrace{\int_{\mathcal{K}_h^i} (a_j a_j d_i u)(d_j a_i a_j L^{ij})(d_i u)}_{\text{1st case}}$$

$$\simeq - \int_{\mathring{\mathcal{K}}_h} (d_j a_j u) d_i[(a_i a_j a_j L^{ij})(d_i u)] \quad (\text{since } a_j d_j u = 0 \text{ on } \partial\mathcal{K}_h)$$

$$= - \int_{\mathring{\mathcal{K}}_h} (d_j a_j u)(a_i a_i a_j a_j L^{ij})(d_i d_i u) - \underbrace{\int_{\mathring{\mathcal{K}}_h} (d_j a_j u)(d_i a_i a_j a_j L^{ij})(d_i u)}_{\text{1st case}}$$

Finally, using (2.4), we have, since $|h^2 dd(L)| \leq C|hd(L)|$

$$|a_i a_i(L) - L| = \left|\frac{h^2}{4} d_i d_i(L)\right| \leq Ch\epsilon_d \quad \text{on } \mathcal{K}_h,$$

so that

$$\int_{\mathring{\mathcal{K}}_h} (d_j a_j u)(a_i a_i a_j a_j L^{ij})(d_i d_i u) \simeq \int_{\mathring{\mathcal{K}}_h} (d_j a_j u) L^{ij}(d_i d_i u).$$

Gathering these equations yields

$$\int_{\mathring{\mathcal{K}}_h} (d_j a_j u) L^{ij}(d_i d_i u) \simeq -\int_{\mathring{\mathcal{K}}_h} (d_j a_j u) L^{ij}(d_i d_i u),$$

up to an error term of the form

$$C\epsilon_d \|u\|_{\mathring{H}^1(\mathcal{B}_h)}^2,$$

which proves the required result thanks to Proposition 3.4. □

*3.4.2. Proof of (3.19) in Lemma 3.11* The idea of the proof is the same as for (3.18), except that the scaling here is $C|s|^{-2}\epsilon_d(\|Su\|_{L^2(\mathcal{K}_h)}^2 + \|Au\|_{L^2(\mathcal{K}_h)}^2)$. Again, this scaling comes from error term corresponding to two cases:

$$h^2 d(L) d(u) d(u) \text{ (first case) or } u d(L) u \text{ (second case)}.$$

Indeed $\|d(L)\|_{L^\infty} \leq C\epsilon_d$, $\|hd(u)\| \leq C\|u\|$ and $\|u\|_{L^2} \leq C|s|^{-1}(\|Su\| + \|Au\|)$.

Below, the symbol $\simeq$ will refer to an equality modulo a term $C|s|^{-2}\epsilon_d(\|Su\|^2 + \|Au\|^2)$.

$$\int_{\mathring{\mathcal{K}}_h} u L^{ij}(a_j d_j u) = \int_{\mathcal{K}_h^j} a_j(u L^{ij})(d_j u) \text{ (since } u = 0 \text{ on } \partial\mathcal{K}_h\text{)}$$

$$= \int_{\mathcal{K}_h^j} (a_j u)(a_j L^{ij})(d_j u) + \underbrace{\frac{h^2}{4} \int_{K^j} (d_j u)(d_j L^{ij})(d_j u)}_{\text{1st case}}$$

$$\simeq -\int_{\mathring{\mathcal{K}}_h} d_j[(a_j u)(a_j L^{ij})] u \text{ (since } u = 0 \text{ on } \partial\mathcal{K}_h\text{)}$$

$$= -\int_{\mathring{\mathcal{K}}_h} (d_j a_j u)(a_j a_j L^{ij}) u - \underbrace{\int_{\mathring{\mathcal{K}}_h} (a_j a_j u)(d_j a_j L^{ij}) u}_{\text{2nd case}}.$$

As in the previous section, using

$$|aa(L) - L| = \left|\frac{h^2}{4} dd(L)\right| \leq Ch\epsilon_d \text{ on } \mathcal{K}_h,$$

we get

$$\int_{\mathring{\mathcal{K}}_h} (d_j a_j u) L^{ij} u \simeq -\int_{\mathring{\mathcal{K}}_h} u L^{ij}(d_j a_j u)$$

which proves the required result. □

*3.5. Proof of Theorem 3.1*

As explained in the discussion following the statement of Theorem 3.1, we base our analysis on the decomposition:

$$\|\Delta_{s,h}u_h\|^2_{L^2(\mathcal{K}_h)} = \|S_{s,h}u_h\|^2_{L^2(\mathcal{K}_h)} + \|A_{s,h}u_h\|^2_{L^2(\mathcal{K}_h)} + 2\int_{\mathcal{K}_h}(S_{s,h}u_h)(A_{s,h}u_h)$$

First set $c > 0$ such that Proposition 3.4 holds for $\max\{h, \epsilon_d, \epsilon_a\} \leq c$ and $|s| \leq c\min\{h^{-2/3}, \epsilon_d^{-1}\}$ with some constant $C$

Then set $\varepsilon > 0$ small enough compared to $C$ so that the error term, estimated in Lemma 3.11 can be absorbed.

Thus, reducing $c$ if necessary, we take $c < c_\varepsilon$ and for $|s| \leq c\min\{\epsilon_d^{-1}, h^{-2/3}\}$, we obtain the Carleman estimate (3.2). □

*3.6. The case $k = d$*

Following Remarks 3.7 and 3.10 and the above proof of Theorem 3.1, one easily proves that in the case $k = d$, Theorem 3.1 holds for $s \in \mathbb{R}^d$ such that $|s| \leq c\min\{\epsilon_d^{-1}, h^{-1}\}$ (and still $(h, \epsilon_d, \epsilon_a) \in (0, c)^3$).

*From now on, we fix $c > 0$ such that Theorem 3.1 holds and we implicitly restrict ourselves to parameters $(h, \epsilon_d, \epsilon_a) \in (0, c)^3$.*

## 4. Construction of the CGO solutions

In the continuous case, CGO solutions, originally introduced in [11] are solutions of the continuous elliptic equation (1.2) that behave like $u(x) \simeq e^{\eta \cdot x}$ for $\eta \in \mathbb{C}^d$, $\eta \cdot \eta = 0$ and $|\eta|$ large enough.

In the continuous case, the existence of a CGO solution relies on two facts :

- The function $u(x) = e^{\eta \cdot x}$ solves $\Delta u = 0$ when $\eta \in \mathbb{C}^d$ satisfies $\eta \cdot \eta = 0$.
- The multiplication by $q$ is a zero order operator that is negligible with respect to the Laplace operator $\Delta$ when $|\eta|$ is big enough. This fact is measured via the Carleman estimate.

This section is an adaptation of the proof of the continuous case, the main difference is that the function $u_h(x) = e^{\eta \cdot x}$ is not solution of $\Delta_h u_h = 0$ anymore for $\eta \in \mathbb{C}^d$ satisfying $\eta \cdot \eta = 0$.

*4.1. CGO solutions are almost harmonic*

Let us begin our analysis by giving estimates on the lack of "discrete harmonicity" of the functions $u(x) = e^{\eta \cdot x}$ for $\eta \in \mathbb{C}^d$ satisfying $\eta \cdot \eta = 0$:

**Proposition 4.1** *For all $C > 0$, there exists $C_0 > 0$ such that for any vector $\eta \in \mathbb{C}^d$ such that $\eta \cdot \eta = 0$ and $|\eta| \leq Ch^{-1}$, define $\Phi_\eta(x) = \eta \cdot x \in C(\mathcal{K}_h)$,*

$$\|e^{-\Phi_\eta}\Delta_h e^{\Phi_\eta}\|_{L^\infty(\check{\mathcal{K}}_h)} \leq C_0\left(|\eta|\epsilon_d + |\eta|^2\epsilon_a + |\eta|^4 h^2\right).$$

**Proof** Using Lemma 2.2, and the fact that $d_i d_i(\eta \cdot x) = 0$, we have

$$e^{-\Phi_\eta} d_i d_i(e^{\Phi_\eta}) = \frac{4}{h^2} sh^2\left(\frac{h\eta \cdot e_i}{2}\right) = (\eta \cdot e_i)^2 + O(|\eta|^4 h^2),$$

$$\left|e^{-\Phi_\eta} a_i d_i(e^{\Phi_\eta})\right| = \left|\frac{2}{h} ch(\eta \cdot e_i) sh\left(\frac{h\eta \cdot e_i}{2}\right)\right| \leq C|\eta|.$$

.

Hence, recalling that $|d_i(\sigma_h^i)| \leq \epsilon_d$,

$$e^{-\Phi_\eta} \sum_i d_i \sigma^i d_i e^{\Phi_\eta} = \sum_i a_i(\sigma^i)(e^{-\Phi_\eta} d_i d_i e^{\Phi_\eta}) + d_i(\sigma^i)(e^{-\Phi_\eta} a_i d_i e^{\Phi_\eta})$$

$$= \sum_i a_i(\sigma^i)(\eta \cdot e_i)^2 + O(|\eta|\epsilon_d + |\eta|^4 h^2).$$

Using (2.8) and $\eta \cdot \eta = 0$, we thus obtain

$$\left|e^{-\Phi_\eta} \sum_i d_i \sigma^i d_i e^{\Phi_\eta}\right| \leq C(|\eta|\epsilon_d + |\eta|^2 \epsilon_a + |\eta|^4 h^2),$$

which concludes the proof of Proposition 4.1. □

*4.2. An auxiliary problem*

In this section, we prove the existence of a solution to an elliptic problem that will arise naturally when proving the existence of discrete CGO solutions.

**Lemma 4.2** *Let $m \in \mathbb{R}_+$, for all $q_h \in C(\mathring{\mathcal{B}}_h)$ satisfying $\|q_h\|_{L^\infty(\mathring{\mathcal{B}}_h)} \leq m$, there exists $s_0 > 0$ depending on $m$ such that if $s_0 \leq |s| \leq c\min\{\epsilon_d^{-1}, h^{-2/3}\}$, then for all $f_h \in C(\mathring{\mathcal{B}}_h)$, there exists a solution $u_h \in C(\mathcal{B}_h)$ to the problem*

$$\Delta_{s,h} u_h + q_h u_h = f_h \quad \text{on } \mathring{\mathcal{B}}_h, \tag{4.1}$$

*(recall that $\Delta_{s,h}$ is defined in (3.1)) that verifies,*

$$\|u_h\|_{L^2(\mathcal{B}_h)} \leq \frac{C}{|s|} \|f_h\|_{L^2(\mathring{\mathcal{B}}_h)}. \tag{4.2}$$

**Proof** Define the operator $P_{s,h}^\star : C_c(\mathcal{B}_h) \to C(\mathcal{B}_h)$ by

$$P_{s,h}^\star v_h := \Delta_{-s,h} v_h + q_h v_h, \quad v_h \in C_c(\mathcal{B}_h).$$

Note that $P_{s,h}^\star$ is the adjoint of the operator $P_{s,h}$ defined by $P_{s,h}(u_h) = \Delta_{s,h} u_h + q_h u_h$.

Let us then consider the following minimization problem: Find $V_h \in C_c(\mathcal{B}_h)$ that minimizes the functional $J_{s,h}$ defined by

$$J_{s,h}(v_h) = \int_{\mathcal{B}_h} |P_{s,h}^\star(v_h)|^2 - \int_{\mathring{\mathcal{B}}_h} f_h v_h$$

among all $v_h \in C_c(\mathcal{B}_h)$.

The existence of such a minimizer is ensured by the fact that $C_c(\mathcal{B}_h)$ is a finite dimensional set and the Carleman estimate in Theorem 3.1 that reads, for $v_h \in C_c(\mathcal{B}_h)$,

$$\|P_{s,h}^\star v_h\|_{L^2(\mathcal{B}_h)}^2 \geq \frac{1}{2} \|\Delta_{-s,h} v_h\|_{L^2(\mathcal{B}_h)}^2 - \|q_h\|_{L^\infty(\mathring{\mathcal{B}}_h)}^2 \|v_h\|_{L^2(\mathcal{B}_h)}^2$$

$$\geq C^{-1} \left(|s|^2 \|v_h\|_{L^2(\mathcal{B}_h)}^2 + \|v\|_{\mathring{H}^1(\mathcal{B}_h)}^2\right) - m^2 \|v_h\|_{L^2(\mathcal{B}_h)}^2$$

$$\geq C^{-1} \left(|s|^2 \|v_h\|_{L^2(\mathcal{B}_h)}^2 + \|v\|_{\mathring{H}^1(\mathcal{B}_h)}^2\right), \tag{4.3}$$

for $|s|$ large enough compared to $m$ and $C$.

As a critical point of $J_{s,h}$, $V_h$ satisfies the following Euler-Lagrange equation:

$$\int_{\mathcal{B}_h} P^\star_{s,h} V_h P^\star_{s,h} v_h = \int_{\mathring{\mathcal{B}}_h} f_h v_h, \quad \forall v_h \in C_c(\mathcal{B}_h). \tag{4.4}$$

Therefore, setting $u_h = P^\star_{s,h} V_h$, doing integration by parts,

$$\int_{\mathcal{B}_h} v_h \left( P_{s,h} u_h - f_h \right) = 0, \quad \forall v_h \in C_c(\mathcal{B}_h),$$

and thus $u_h$ solves (4.1). Here, we have used that, since $v_h \in C_c(\mathcal{B}_h)$ and $u_h \in C_c(\mathcal{K}_h)$ since $\mathcal{B}_h \subset \mathring{\mathcal{K}}_h$,

$$\int_{\mathcal{B}_h} u_h P^\star_{s,h} v_h = \int_{\mathcal{K}_h} u_h P^\star_{s,h} v_h = \int_{\mathcal{K}_h} P_{s,h} u_h v_h = \int_{\mathcal{B}_h} P_{s,h} u_h v_h.$$

To get the estimate on the $L^2(\mathcal{B}_h)$-norm of $u_h$, we plug $v_h = V_h$ in (4.4) and use (4.3):

$$\|P^\star_{s,h} V_h\|^2_{L^2(\mathcal{B}_h)} = \int_{\mathring{\mathcal{B}}_h} f_h V_h \leq \|f_h\|_{L^2(\mathring{\mathcal{B}}_h)} \|V_h\|_{L^2(\mathcal{B}_h)}$$
$$\leq \frac{C}{|s|} \|f_h\|_{L^2(\mathring{\mathcal{B}}_h)} \|P^\star_{s,h} V\|_{L^2(\mathcal{B}_h)}.$$

We then have

$$\|u_h\|_{L^2(\mathcal{B}_h)} = \|P^\star_{-s,h} V_h\|_{L^2(\mathcal{B}_h)} \leq \frac{C}{|s|} \|f_h\|_{L^2(\mathring{\mathcal{B}}_h)},$$

which is exactly the sought result. □

We shall also need further estimate on the function $u_h$ given by Lemma 4.2, namely in the $H^1$-norm.

At this stage, let us recall Assumption 4 that states the existence of a suitable cut-off function $\psi$ that equals to 1 on $\mathcal{W}_h$ while being compactly supported in $\mathcal{B}_h$ with additional estimates on its discrete derivatives.

**Lemma 4.3** *Let $m \in \mathbb{R}_+$, for all $q_h \in C(\mathring{\mathcal{W}}_h)$ satisfying $\|q_h\|_{L^\infty(\mathring{\mathcal{W}}_h)} \leq m$, there exists $s_0 > 0$ depending on $m$ such that if $s_0 \leq |s| \leq c_0 \min\{\epsilon_d^{-1}, h^{-2/3}\}$, for all $f_h \in C(\mathring{\mathcal{W}}_h)$, there exists a solution $u_h \in C(\mathcal{W}_h)$ to the problem*

$$\Delta_{s,h} u_h + q_h u_h = f_h \quad \text{on } \mathring{\mathcal{W}}_h, \tag{4.5}$$

*that verifies,*

$$\|u_h\|_{\dot{H}^1(\mathcal{W}_h)} + |s| \|u_h\|_{L^2(\mathcal{W}_h)} \leq C \|f_h\|_{L^2(\mathring{\mathcal{W}}_h)}. \tag{4.6}$$

**Proof** First, we extend $f_h$ and $q_h$ to $\mathcal{B}_h$ by zero outside $\mathring{\mathcal{W}}_h$, and we define $u_h$ as in Lemma 4.2.

Thus, to prove Lemma 4.3, it remains to show that $\|u_h\|_{\dot{H}^1(\mathcal{W}_h)} \leq C \|f_h\|_{L^2(\mathring{\mathcal{W}}_h)}$ since we already have by Lemma 4.2 that

$$|s| \|u_h\|_{L^2(\mathcal{W}_h)} \leq |s| \|u_h\|_{L^2(\mathcal{B}_h)} \leq C \|f\|_{L^2(\mathring{\mathcal{W}}_h)}.$$

Thus, to get the estimate in the $\dot{H}^1(\mathcal{W}_h)$-norm, we do a multiplier-type estimate, multiplying the equation satisfied by $u_h$ by $\psi u_h$, where $\psi$ is the multiplier function given by Assumption 4.

First, let us remark that the $L^2(\mathring{\mathcal{B}}_h)$-norm of $u_h$ is bounded by $C\|f_h\|_{L^2(\mathring{\mathcal{W}}_h)}/|s|$ and the one of $\Delta_{s,h} u = f - qu$ is then bounded by $C\|f\|_{L^2(\mathring{\mathcal{W}}_h)}$, and then

$$\left| \int_{\mathring{\mathcal{B}}_h} \Delta_{s,h} u_h \psi u_h \right| \leq \frac{C}{|s|} \|f_h\|_{L^2(\mathring{\mathcal{W}}_h)}. \tag{4.7}$$

But on the other hand, following the computation already performed in the proof of the Proposition 3.3, we have:

$$\int_{\mathring{\mathcal{B}}_h} \Delta_{s,h} u_h \psi u_h = \sum_i \int_{\mathcal{B}_h^i} \left( r_1^i d_i(u_h) d_i(\psi u_h) + r_2^i a_i(u_h) a_i(\psi u_h) \right.$$
$$\left. + r_3^i d_i(u_h) a_i(\psi u_h) - r_3^i a_i(u_h) d_i(\psi u_h) \right), \tag{4.8}$$

where the coefficients $r_1^i, r_2^i, r_3^i$ are given in (3.6)–(3.7)–(3.8):

$$r_1^i = \sigma^i ch^2\left(\frac{h}{2} s \cdot e_i\right) = \sigma^i + O(h^2 |s|^2), \tag{4.9}$$

$$r_2^i = \frac{-4\sigma^i}{h^2} sh^2\left(\frac{h}{2} s \cdot e_i\right), \tag{4.10}$$

$$r_3^i = \frac{-\sigma^i}{h} sh(hs \cdot e_i). \tag{4.11}$$

In particular, for $|s|h$ bounded,

$$|r_1^i| \leq C, \quad |r_2^i| \leq C|s|^2, \quad |r_3^i| \leq C|s|. \tag{4.12}$$

for some constant $C$ independent of $s$ and $h$.

Combining (4.7) and (4.8) and using the formula of lemma 2.1, we obtain

$$\frac{C}{|s|} \|f_h\|_{L^2(\mathring{\mathcal{W}}_h)}^2 \geq \sum_i \int_{\mathcal{B}_h^i} \underbrace{(d_i u_h)^2 \left( r_1^i a_i(\psi) + r_3^i \frac{h^2}{4} d_i(\psi) \right)}_{(1)}$$
$$+ \sum_i \int_{\mathcal{B}_h^i} \underbrace{d_i(u_h) a_i(u_h) \left( r_1^i - r_2^i \frac{h^2}{4} \right) d_i(\psi)}_{(2)}$$
$$+ \sum_i \int_{\mathcal{B}_h^i} \underbrace{(a_i(u_h))^2 \left( r_2^i a_i(\psi) - r_3^i d_i(\psi) \right)}_{(3)}. \tag{4.13}$$

We shall now estimate the terms **(1)**, **(2)** and **(3)** separately.

To estimate **(1)**, we use the fact that $|d_i(\psi)| \leq C$ and $|r_3^i| \leq C|s|$ by (4.11), so that $|r_3^i h^2 d_i(\psi)| \leq Ch^2|s|$. Now, we use that the operator $d_i$ is of norm of the order of $1/h$, hence

$$\sum_{i=1}^k \int_{\mathcal{B}_h^i} (d_i u_h)^2 \leq \frac{C}{h^2} \int_{\mathcal{B}_h} (u_h)^2. \tag{4.14}$$

Combining all these estimates, we obtain

$$\left| \sum_i \int_{\mathcal{B}_h^i} (d_i u_h)^2 r_3^i \frac{h^2}{4} d_i(\psi) \right| \leq C \|f_h\|_{L^2(\mathring{\mathcal{W}}_h)}^2 |s|^{-1}.$$

According to (4.9), we also have

$$\left| \sum_i \int_{\mathcal{B}_h^i} (d_i u_h)^2 \, r_1^i a_i(\psi) - \sum_i \int_{\mathcal{B}_h^i} (d_i u_h)^2 \, \sigma_h^i a_i(\psi) \right|$$
$$\leq C|s|^2 |h|^2 \|u_h\|_{\dot{H}^1(\mathcal{B}_h)}^2 \leq C|s|^2 \|u_h\|_{L^2(\mathcal{B}_h)}^2 \leq C\|f_h\|_{L^2(\mathring{\mathcal{W}}_h)}^2.$$

where we again used (4.14).

Finally, using that $a_i(\psi) \geq 0$ on $\mathcal{K}_h$ and $a_i(\psi) = 1$ on $\mathcal{W}_h$, we obtain

$$\sum_i \int_{\mathcal{B}_h^i} (d_i(u_h))^2 \left( r_1^i a_i(\psi) + r_3^i \frac{h^2}{4} d_i(\psi) \right) \geq \|u_h\|_{\dot{H}^1(\mathcal{W}_h)}^2 - C \frac{\|f_h\|_{L^2(\mathring{\mathcal{W}}_s)}^2}{|s|}, \quad (4.15)$$

The term denoted as (**2**) simplifies by remarking that $a_i(u)d_i(u) = d_i(u^2/2)$ and by performing an integration by part without boundary terms since $d_i\psi = 0$ on $\partial \mathcal{B}_h$. To simplify notations, let us also remark that $r_1^i - r_2^i h^2/4 = \sigma_h^i \, ch(hs \cdot e_i)$. Then

$$\sum_i \int_{\mathcal{B}_h^i} d_i(u_h) a_i(u_h) \sigma_h^i ch(hs \cdot e_i) d_i(\psi) = - \sum_i \int_{\mathcal{B}_h^i} \frac{u_h^2}{2} d_i \left( \sigma_h^i ch(hs \cdot e_i) d_i(\psi) \right).$$

Recalling that $|d_i(\sigma^i d_i \psi)| \leq C$ by assumption, this term then satisfies

$$\left| \sum_i \int_{\mathcal{B}_h^i} d_i(u_h) a_i(u_h) \sigma_h^i ch(hs \cdot e_i) d_i(\psi) \right| \leq C\|u_h\|_{L^2(\mathcal{B}_h)}^2 \leq C \frac{\|f_h\|_{L^2(\mathring{\mathcal{W}}_s)}^2}{|s|^2}. \quad (4.16)$$

Finally, to estimate the term (**3**), using the rough estimates (4.12), one easily checks that

$$\left| \sum_i \int_{\mathcal{B}_h^i} (a_i(u_h))^2 \left( r_2^i a_i(\psi) - r_3^i d_i(\psi) \right) \right| \leq C|s|^2 \|u_h\|_{L^2(\mathcal{B}_h)}^2 \leq C\|f_h\|_{L^2(\mathring{\mathcal{W}}_s)}^2. \quad (4.17)$$

Combining estimates (4.15)–(4.16)–(4.17) and putting them in (4.13), we obtain directly estimate (4.6). □

### 4.3. Existence of discrete CGO solutions

In this section we state the theorem of existence of discrete CGO solutions and we give some estimates on them.

**Theorem 4.4** *Let $m \in \mathbb{R}_+$. For all $q_h \in C(\mathring{\mathcal{W}}_h)$ satisfying $\|q_h\|_{L^\infty(\mathring{\mathcal{W}}_h)} \leq m$, there exist $s_0 > 0$ that depends on $m$ such that $\forall \eta \in \mathbb{C}^d$ such that $\eta \cdot \eta = 0$, if $s := \Re(\eta)$ verifies $s_0 \leq |s| \leq c \min\{\epsilon_d^{-1}, h^{-2/3}\}$, there exists $u_h \in C(\mathcal{W}_h)$ a solution of*

$$\Delta_h u_h + q_h u_h = 0, \quad \text{on } \mathring{\mathcal{W}}_h,$$

*that satisfies*

$$u_h(x) = e^{\eta \cdot x} + e^{s \cdot x} r_h(x) \text{ on } \mathcal{W}_h \quad (4.18)$$

*with*

$$\|r_h\|_{\dot{H}^1(\mathcal{W}_h)} + |s|\|r_h\|_{L^2(\mathcal{W}_h)} \leq C(1 + |s|^2 \epsilon_a + |s|^4 h^2), \quad (4.19)$$

*Moreover, the CGO solution $u_h$ satisfies,*

$$\|u_h\|_{H^1(\mathcal{W}_h)} \leq C e^{|s|} |s|^2. \quad (4.20)$$

Before going into the proof, let us remark that when $\eta \in \mathbb{C}^d$ satisfies $\eta \cdot \eta = 0$, easy computations show $|\Re(\eta)| = |\Im(\eta)|$ and then $|s| = |\eta|/\sqrt{2}$.

**Proof** Set $\Phi_\eta(x) = \eta \cdot x \in C(\mathcal{K}_h)$ and $\phi_s(x) = s \cdot x = \Re(\Phi_\eta(x))$.

Writing $\eta = s + i\alpha$, $(s, \alpha) \in (\mathbb{R}^d)^2$, using Proposition 4.1, we get

$$e^{-\phi_s}(\Delta_h e^{\Phi_\eta} + q_h e^{\Phi_\eta}) = e^{i\alpha \cdot x}\left(e^{-\Phi_\eta}\Delta e^{\Phi_\eta} + q\right)$$

and therefore

$$\|e^{-\phi_s}(\Delta_h e^{\Phi_\eta} + q_h)\|_{L^2(\mathring{\mathcal{W}}_h)} \leq C(|s|^2\epsilon_a + |s|^4 h^2 + 1),$$

where, since $|s|\epsilon_d \leq c_0$, we have bounded $|s|\epsilon_d$ by some uniform constant.

The sought $r_h \in C(\mathcal{W}_h)$ must verify

$$\Delta_h\left(e^{\phi_s}r\right) + q_h e^{\phi_s}r_h = -\Delta_h e^{\Phi_\eta} - q_h e^{\Phi_\eta} \quad \text{in } \mathcal{W}_h,$$

that is

$$\Delta_{s,h}r_h + q_h r_h = e^{-\phi_s}(\Delta_h e^{\Phi_\eta} + q_h e^{\Phi_\eta}) \quad \text{in } \mathcal{W}_h.$$

Lemmas 4.2–4.3 then prove the existence of such an $r_h$ that verifies in addition

$$\|r_h\|_{\dot{H}^1(\mathcal{W}_h)} + |s|\|r_h\|_{L^2(\mathcal{W}_h)} \leq \|e^{-\phi_s}(\Delta_h e^{\Phi_\eta} + q_h)\|_{L^2(\mathring{\mathcal{W}}_h)}$$
$$\leq C(1 + |s|^2\epsilon_a + |s|^4 h^2),$$

which ends the proof of (4.19).

Estimate (4.20) then immediately follows. □

Of course, replacing $r_h$ by $\tilde{r}_h = e^{-i\Im(\eta) \cdot x}r_h$, the solution $u_h$ in (4.18) can be rewritten as

$$u_h(x) = e^{\eta \cdot x}(1 + \tilde{r}_h(x)) \text{ on } \mathcal{W}_h, \tag{4.21}$$

and easy computations show that $\tilde{r}_h$ also satisfy estimates (4.19).

In the sequel, we shall use this form of CGO solutions instead of (4.18).

**Remark 4.5** *Following Section 3.6, when $k = d$, one can construct CGO solutions for $|s| \leq c\min\{\epsilon_d^{-1}, h^{-1}\}$. Indeed, Lemma 4.2 and 4.3 only require $|s| \leq ch^{-2/3}$ in order to apply Theorem 3.1, which holds when $k = d$ up to the scaling $h^{-1}$.*

*However, in this case, we see that $\|r_h\|_{L^2}$ is bounded by $|s|^3 h^2$, which is small only for $|s|h^{2/3}$ small, thus explaining why we cannot improve much the uniform stability results when $k = d$ despite the fact that the Carleman estimate is better in that case.*

*Howver, we shall see in Section 6 that, in the case of a uniform mesh, using special "discrete harmonic functions", we can slightly improve the stability estimates.*

### 5. Stability estimate

In this Section, we focus on the proof of Theorem 2.7.

**Proof** Let $q_{1,h}$ and $q_{2,h}$ be two potentials in $C(\mathring{\mathcal{W}}_h)$ such that $\|q_{j,h}\|_{L^\infty(\mathring{\mathcal{W}}_h)} \leq m$, $j = 1, 2$.

For any $\xi \in \hat{\mathcal{K}}_h$ satisfying $2\pi|\xi| \leq c\min\{\epsilon_d^{-1}, h^{-2/3}\}$, where $c > 0$ is the one of Theorem 4.4, we set $\beta = -2\pi\xi \in \mathbb{R}^d$.

Then find $s$ and $\delta$ in $\mathbb{R}^d$ such that

- $s_0 \leq |s| \leq c\min\{\epsilon_d^{-1}, h^{-2/3}\}$, where $s_0$ is given by Theorem 4.4;
- the three real vectors $(s, \beta, \delta)$ are orthogonal: $s \cdot \beta = \beta \cdot \delta = \delta \cdot s = 0$;
- $|s|^2 = |\delta|^2 + |\beta|^2$.

For this to be possible, we strongly use the fact that we are in dimension greater than 3, $d \geq 3$. Also note that these conditions can be satisfied if and only if $|s| \geq |\beta|$.

Set then $\eta_1$ and $\eta_2$ as

$$\eta_1 = s + i(\beta + \delta) \text{ and } \eta_2 = -s + i(\beta - \delta).$$

By construction, $\eta_1$ and $\eta_2$ obeys the condition of Theorem 4.4. Therefore, there exist solutions $u_{1,h}$ and $u_{2,h}$ of, respectively,

$$\Delta_h u_{j,h} + q_{j,h} u_{j,h} = 0, \quad \text{on } \mathring{\mathcal{W}}_h, \ j = 1, 2,$$

that verify

$$u_j(x) = e^{\eta_j \cdot x}(1 + r_j), \quad j = 1, 2.$$

In order to simplify the notations, set $\lambda := \|\Lambda_h[q_{1,h}] - \Lambda_h[q_{2,h}]\|_{\mathfrak{L}_h}$.

Equation (2.24) then reads

$$\int_{\mathring{\mathcal{W}}_h} (q_{1,h} - q_{2,h}) e^{i\beta \cdot x}(1 + r_{1,h})(1 + r_{2,h}) = \int_{\mathring{\mathcal{W}}_h} (q_{1,h} - q_{2,h}) u_{1,h} u_{2,h}$$

$$= \int_{\partial \mathcal{W}_h} (\Lambda_h[q_{1,h}] - \Lambda_h[q_{2,h}])(u_{1,h}) u_{2,h}$$

$$\leq \lambda |u_{1,h}|_{H^{1/2}(\partial \mathcal{W}_h)} |u_{2,h}|_{H^{1/2}(\partial \mathcal{W}_h)}$$

$$\leq \lambda \|u_{1,h}\|_{H^1(\mathcal{W}_h)} \|u_2\|_{H^1(\mathcal{W}_h)} \leq \lambda |s|^2 e^{2|s|},$$

where we used (4.20) in the last inequality. Therefore,

$$|\mathcal{F}_h(q_{1,h} - q_{2,h})(\xi)| \leq \lambda s^2 e^{2s} + \int_{\mathring{\mathcal{W}}_h} (q_1 - q_2) e^{i\beta \cdot x}(r_{1,h} + r_{2,h} + r_{1,h} r_{2,h})$$

$$\leq \lambda s^2 e^{2s} + C(\|r_{1,h}\|_{L^2(\mathcal{W}_h)} + \|r_{2,h}\|_{L^2(\mathcal{W}_h)} + \|r_{1,h}\|_{L^2(\mathcal{W}_h)} \|r_{2,h}\|_{L^2(\mathcal{W}_h)}) \quad (5.1)$$

We first suppose that $s$ is chosen so that it also verifies $|s| \leq \epsilon_a^{-1}$. This implies that $\|r_{1,h}\|_{L^2(\mathcal{W}_h)}$ and $\|r_{2,h}\|_{L^2(\mathcal{W}_h)}$ are bounded: indeed, using (4.19),

$$\|r_{1,h}\|_{L^2(\mathcal{W}_h)} + \|r_{2,h}\|_{L^2(\mathcal{W}_h)} \leq C \left( \frac{1}{|s|} + |s|^3 h^2 + \epsilon_a |s| \right).$$

Setting $\tilde{\mu} = \max\{\epsilon_a^{1/2}, h^{1/2}, \epsilon_d/c\}$, we always have $1/\tilde{\mu} \leq c \min\{\epsilon_d^{-1}, h^{-2/3}\}$ and a simple study shows that, setting $|s| = 1/\tilde{\mu}$,

$$\left( \frac{1}{|s|} + |s|^3 h^2 + \epsilon_a |s| \right) \leq C\tilde{\mu}.$$

Therefore, if $\lambda = 0$, for all $\xi \in \hat{\mathcal{K}}_h$ such that $2\pi|\xi| \leq 1/\tilde{\mu}$,

$$|\mathcal{F}_h(q_{1,h} - q_{2,h})(\xi)| \leq C\tilde{\mu}.$$

Of course, using (5.1), this is still the case when considering $\lambda$ small enough, namely such that

$$\frac{\lambda}{\tilde{\mu}^2} e^{2/\tilde{\mu}} \leq \tilde{\mu},$$

which can be guaranteed, for instance, for

$$-\log(\lambda) \geq \frac{3}{\tilde{\mu}}. \tag{5.2}$$

Therefore, if (5.2) holds, taking $|s| = 1/\tilde{\mu}$, we obtain that for all $\xi \in \hat{\mathcal{K}}_h$ such that $2\pi|\xi| \leq 1/\tilde{\mu}$,

$$|\mathcal{F}_h(q_{1,h} - q_{2,h})(\xi)| \leq C\tilde{\mu}. \tag{5.3}$$

On the other hand, if $-\log(\lambda) \in (4|s_0|, 3/\tilde{\mu})$, taking $|s| = -\log(\lambda)/3$ in (5.1), we get, for all $\xi \in \hat{\mathcal{K}}_h$ such that $2\pi|\xi| \leq |\log(\lambda)|/3$,

$$|\mathcal{F}_h(q_{1,h} - q_{2,h})(\xi)| \leq C\lambda^{1/3}(\log(\lambda))^2 + \frac{1}{|\log(\lambda)|} + |\log(\lambda)|^3 h^2 + |\log(\lambda)|\epsilon_a$$
$$\leq \frac{C}{|\log(\lambda)|},$$

where we used that $\tilde{\mu} = \max\{h^{1/2}, \epsilon_a^{1/2}, \epsilon_d/c_0\} \leq 3/|\log(\lambda)|$,

Combining these two cases, we obtain that, if $\lambda \leq e^{-4s_0}$, setting

$$\mu = \max\{\tilde{\mu}, 3/\log(\lambda)\},$$

for all $\xi \in \hat{\mathcal{K}}_h$ such that $2\pi|\xi| \leq 1/\mu$,

$$|\mathcal{F}_h(q_{1,h} - q_{2,h})(\xi)| \leq C\mu, \tag{5.4}$$

which of course coincides with (2.21).

We now turn our attention to proving the bounds (2.22) on the $H^{-r}$ norms of $q_{1,h} - q_{2,h}$.

For any $\rho \in (0, 1/\mu)$, we have

$$|q_{1,h} - q_{2,h}|^2_{H^{-r}(\mathcal{K}_h)} = \sum_{\xi \in \hat{\mathcal{K}}_h} |\mathcal{F}_h(q_{1,h} - q_{2,h})(\xi)|^2 (1 + |\xi|^2)^{-r}$$
$$= \sum_{|\xi| < \rho, \xi \in \hat{\mathcal{K}}_h} |\mathcal{F}_h(q_{1,h} - q_{2,h})(\xi)|^2 (1 + |\xi|^2)^{-r}$$
$$+ \sum_{|\xi| > \rho, \xi \in \hat{\mathcal{K}}_h} |\mathcal{F}_h(q_{1,h} - q_{2,h})(\xi)|^2 (1 + |\xi|^2)^{-r}$$
$$\leq C(\rho^d \mu^2 + \rho^{-2r}),$$

where we used that the two potentials $q_{1,h}$ and $q_{2,h}$ have $L^2$-norms bounded by a constant $m$ and therefore,

$$\sum_{\xi \in \hat{\mathcal{K}}_h} |\mathcal{F}_h(q_{1,h} - q_{2,h})(\xi)|^2 \leq C.$$

Optimizing in $\rho \leq 1/\mu$, we take $\rho = \mu^{-2/(d+2r)}$, which is indeed smaller than $1/\mu$, and we obtain (2.22). □

**Remark 5.1** *Let us emphasize that the above proof requires the knowledge of the norm of the difference between the DtN maps only for CGO solutions corresponding to frequency scaling smaller than $\mu^{-1}$.*

*Therefore, we expect Theorem 2.7 to hold when replacing*

$$\|\Lambda_h[q_{1,h}] - \Lambda_h[q_{2,h}]\|_{\mathfrak{L}_h} := \max_{\|g_h\|_{H^{1/2}(\partial \mathcal{W}_h)} = 1} |(\Lambda_h[q_{1,h}] - \Lambda_h[q_{2,h}])(g_h)|_{H^{-1/2}(\partial \mathcal{W}_h)}.$$

*by*

$$\max_{\|u_h\|_{H^1(\mathcal{W}_h)} \leq 1, \, u_h \in V_{\leq \mu^{-1}}} |(\Lambda_h[q_{1,h}] - \Lambda_h[q_{2,h}])((u_h)_{|\partial \mathcal{W}_h})|_{H^{-1/2}(\partial \mathcal{W}_h)},$$

where $V_{\leq \mu^{-1}}$ is the space of functions of frequency smaller than $1/\mu$, defined by
$$V_{\leq \mu^{-1}} = \{u \in C(\mathcal{W}_h) \text{ s.t. } \|u\|_{H^1(\mathcal{W}_h)} \leq \mu^{-1}\|u\|_{L^2(\mathcal{W}_h)}\}.$$

*Unfortunately, we do not know so far if the CGO solutions corresponding to $|\eta| \leq \mu^{-1}$ belong to that space or not. This would be interesting since it would mean that only the solutions that are relevant from the numerical analysis point of view should be taken into account.*

## 6. The special case of uniform meshes

When the meshes are uniform, that is when the set of directions $(e_i)_{i=1..k}$ correspond to an orthonormal basis of $\mathbb{R}^d$ (hence $k = d$) and when $\sigma^i = 1$ for every $i = 1..d$, the analysis greatly simplifies. In this case, we have in particular that $\epsilon_d = \epsilon_a = 0$.

But much more simplifications occur. This section aims at giving results corresponding to that case. In particular, we shall see that in this case, a uniqueness result can be derived from the knowledge of the DtN map, see Theorem 6.5.

### 6.1. Discrete CGO solutions

In the case of a uniform mesh, Theorem 3.1 can be improved as follows:

**Theorem 6.1** *Using the same notations as in Theorem 3.1, if $h \leq c$, for all $|s|$, for all $u_h \in C_c(\mathcal{B}_h)$,*
$$|s|\|u_h\|_{L^2(\mathcal{B}_h)} \leq C\|\Delta_{s,h}u\|_{L^2(\mathcal{K}_h)}. \tag{6.1}$$

**Proof** Of course, the proof closely follows the one of Theorem 3.1. But, in the case of uniform meshes, the proof turns out to be much easier since, using the same decomposition as the one in Proposition 3.3, the commutator $(A_{s,h}u_h, S_{s,h}u_h) = 0$. This fact is due to $d(\sigma) = 0$ implying $d(L) = 0$ as well. We omit the proof since it is based on the same ingredients as in Section 3.4. In order to conclude (6.1), we only need to prove a Poincaré inequality on $A_{s,h}$ that holds for all $s$, which coincides with (3.15) since in this case $\tilde{A}_{s,h} = A_{s,h}$. Details are left to the reader. □

In the uniform case, we can construct more CGO solutions than in the non-uniform case. This is in particular due to the fact that the Carleman estimate (6.1) holds without any limitation on the size of $s$.

In particular, if $\eta \in \mathbb{C}^d$, setting $\Phi_\eta(x) = \eta \cdot x$, one easily checks that $e^{\Phi_\eta}$ is a solution of $\Delta_h u_h = 0$ if and only if
$$\frac{4}{h^2} \sum_{i=1}^{d} sh^2\left(\frac{h\eta \cdot e_i}{2}\right) = 0. \tag{6.2}$$

Thus, based on these discrete harmonic solutions and the Carleman estimate (6.1), we get the following:

**Theorem 6.2** *Let $m \in \mathbb{R}_+$. If $h < c$, for all $q_h \in C(\mathring{\mathcal{W}}_h)$ satisfying $\|q_h\|_{L^\infty(\mathring{\mathcal{W}}_h)} \leq m$, there exist $s_0 > 0$ depending on $m$ such that $\forall \eta \in \mathbb{C}^d$ satisfying (6.2), if $s := \Re(\eta)$ verifies $s_0 < |s|$, there exists $u_h \in C(\mathcal{W}_h)$ a solution of*
$$\Delta_h u_h + q_h u_h = 0, \quad \text{on } \mathring{\mathcal{W}}_h,$$
*that satisfies*
$$u_h(x) = e^{\eta \cdot x}(1 + r_h(x)) \text{ on } \mathcal{W}_h \tag{6.3}$$

*with*

$$|s|\|r_h\|_{L^2(\mathcal{W}_h)} \leq C, \tag{6.4}$$

*for some $C > 0$ independent of $h > 0$ and $s$.*

The proof of Theorem 6.2 is the same as the one of Theorem 4.4, and is therefore omitted.

*6.2. Application*

The application we have in mind is very much related to the proof of Theorem 2.7. Namely, we are going to prove the following:

**Proposition 6.3** *Suppose that the mesh is uniform, that is the set of connection $(e_i)_{i=1..d}$ is an orthonormal set of $\mathbb{R}^d$ and $\sigma^i = 1$ forall $i = 1..d$. Let $m \in \mathbb{R}_+$. Then for all $q_h \in C(\mathring{\mathcal{W}}_h)$ satisfying $\|q_h\|_{L^\infty(\mathring{\mathcal{W}}_h)} \leq m$, for all $\beta \in \mathbb{R}^d$, if there exists $j_0 \in [\![1,d]\!]$ such that $\beta \cdot e_{j_0} = 0$, then for all $h > 0$ there exists $\eta_h \in \mathbb{C}^d$ that can be made arbitrarily large so that there exist $u_{\pm,h} \in C(\mathcal{W}_h)$ solutions of*

$$\Delta_h u_{\pm,h} + q_h u_{\pm,h} = 0,$$

*that satisfy*

$$u_{\pm,h}(x) = e^{i\beta \cdot x} e^{\pm \eta_h \cdot x}(1 + r_{\pm,h}(x))$$

*with $r_{\pm,h}(x)$*

$$|\eta_h|\|r_{\pm,h}\|_{L^2(\mathcal{W}_h)} \leq C.$$

*Furthermore, $\eta_h$ is independent of $q_h$.*

**Proof** Let $\beta \in \mathbb{R}^d$, and assume, without loss of generality, that $\beta \cdot e_1 = 0$. Then, one has to find $\eta_h \in \mathbb{C}^d$ such that $\eta_1 = i\beta + \eta_h$ and $\eta_2 = i\beta - \eta_h$ both satisfy (6.2), that is:

$$0 = \sum_{j=1}^{d} sh^2 \left( \frac{h\eta_h \cdot e_j \pm ih\beta \cdot e_j}{2} \right).$$

After some tedious computations, this yields:

$$0 = \sum_{j=1}^{d} sh(h\eta_h \cdot e_j) \sin(h\beta \cdot e_j) \text{ and } \sum_{j=1}^{d} ch(h\eta_h \cdot e_j) \cos(h\beta \cdot e_j) = d. \tag{6.5}$$

We then choose a real unit vector $a_h$ orthogonal to $\sum_{j=1}^{d} \sin(h\beta \cdot e_j) e_j$ such that $a \cdot e_1 = 0$. Such a vector always exist in dimension $d \geq 3$.

For $\alpha \in \mathbb{R}^*$ arbitrary, let $\eta_h \in \mathbb{C}^d$ be a solution of

$$\forall j = 2..d, \quad \frac{1}{h} sh(h\eta_h \cdot e_j) = \alpha a_{j,h} \tag{6.6}$$

$$\text{and } ch(h\eta_h \cdot e_1) = d - \sum_{j=2}^{d} ch(h\eta_h \cdot e_j) \cos(h\beta \cdot e_j),$$

which can be solved in $\mathbb{C}$ since $ch$ is surjective on $\mathbb{C}$.

By construction, this vector $\eta_h$ is convenient for Proposition 6.3.

Besides the coefficient $\alpha$ may be chosen arbitrarily large enough, which completes the proof of Proposition 6.3. □

**Remark 6.4** *When there is no $j$ such that $\beta \cdot e_j = 0$, we do not know how to solve (6.5) with arbitrarily large $\eta$. This would have important consequences with respect to the uniqueness properties of the discrete Calderón problems on uniform meshes, see for instance the paragraph below.*

*6.3. A uniqueness result*

The above construction allows us to state a uniqueness result:

**Theorem 6.5** *For $q_h$ defined on $\mathcal{K}_h$ and $j = 1..d$, define $q_h[j]$, the average of $q_h$ in the direction $j$ as follows: $q_h[j] : \mathcal{K}_h^{\ominus j} \to \mathbb{R}$, where*

$$\mathcal{K}_h^{\ominus j} := \left\{ x \in [0,1]^{d-1} \text{ such that } \exists k \in [\![0, N-1]\!]^{d-1} \text{ such that } x = \frac{k}{N} \right\},$$

*and is defined by*

$$q_h[j](x_1, .., x_{j-1}, x_{j+1}, .., x_d) :=$$
$$h \sum_{(x_1, x_{j-1}, x_j, x_{j+1}, .., x_d) \in \mathcal{K}_h} q_h(x_1, .., x_{j-1}, x_j, x_{j+1}, .., x_d).$$

*If the mesh is uniform and $\Lambda_h[q_{1,h}] = \Lambda_h[q_{2,h}]$ and $q_{1,h}, q_{2,h}$ belong to $L^\infty(\mathring{\mathcal{W}}_h)$, then for every direction $j$, we have $q_{1,h}[j] = q_{2,h}[j]$.*

**Proof** Fix $j \in [\![1, d]\!]$. Take then $\hat{\xi} \in \hat{\mathcal{K}}_h^{\ominus j}$. Our goal is to show that

$$\mathcal{F}_h(q_{1,h}[j] - q_{2,h}[j])(\hat{\xi}) = 0,$$

which of course implies that $q_{1,h}[j] = q_{2,h}[j]$. In order to do that, we define $\xi \in \hat{\mathcal{K}}_h$ as follows:

$$\xi = (\hat{\xi}_1, \ldots, \hat{\xi}_{j-1}, 0, \hat{\xi}_{j+1}, \ldots, \hat{\xi}_d),$$

One then easily checks that

$$\mathcal{F}_h(q_{1,h} - q_{2,h})(\xi) = \mathcal{F}_h(q_{1,h}[j] - q_{2,h}[j])(\hat{\xi}).$$

Thus we only have to prove that for all $\xi \in \hat{\mathcal{K}}_h$ such that $\xi \cdot e_j = 0$,

$$\mathcal{F}_h(q_{1,h} - q_{2,h})(\xi) = 0.$$

But, setting $\beta = 2\pi\xi$, we can use then the functions $u_+$ and $u_-$ built as in Proposition 6.3 and apply equation (2.24):

$$0 = \int_{\mathcal{W}_h} (q_{1,h} - q_{2,h}) u_{+,h} u_{-,h}$$
$$= \int_{\mathcal{W}_h} (q_{1,h} - q_{2,h}) e^{2i\beta \cdot x}(1 + r_{+,h})(1 + r_{-,h})$$
$$= \mathcal{F}_h(q_{1,h} - q_{2,h})(\xi) + O\left(\frac{1}{|\eta|}\right).$$

Since $\eta$ can be made arbitrarily large, we obtain $\mathcal{F}_h(q_{1,h} - q_{2,h})(\xi) = 0$ for all $\xi \in \hat{\mathcal{K}}_h$ such that $\xi \cdot e_j = 0$. This concludes the proof. □

**7. Conclusion**

In this article, we have derived uniform stability estimates for the discrete Calderón problems. But still, a lot remains to be done.

1. **Convergence of the inverse problems.** The results developed here should be considered as a first step of the convergence of the discrete Calderón problems towards the continuous ones. Indeed, it would be very interesting to prove that if $q_h$ is such that $\Lambda_h[q_h]$ is close to $\Lambda[q]$ for $h$ small enough, then $q_h$ is close to $q$.

However, this might be much more difficult than expected since the operators $\Lambda_h[q_h]$, $\Lambda[q]$ are operators defined on different functional spaces. Besides, as we have explained in that article, $\Lambda_h[q_h]$ contains all the solutions of the discrete Calderón problem, which may have very unexpected behavior at large frequency scales.

2. **Discrete vs continuous CGO.** Another open problem concerns the convergence of the discrete CGO solutions toward the continuous ones. This issue is probably related to the above point, but so far, it is not clear how fast the discrete CGO solutions converge to the continuous ones.

Of course, this study should help to understand in which sense the discrete Calderón problems converge to the continuous one.

3. **The 2 dimensional case.** When considering the Calderón problem in 2-d, one cannot use CGO solutions anymore, and the analysis in [2] rather uses tools of complex analysis. To our knowledge, trying to get uniform stability estimates for the discrete 2-d Calderón problems is completely open.

4. **New bounds on the gradient.** We point out that we obtained new bounds on the gradient of $r$ by performing a multiplier type estimate in the proof of Lemma 4.3. Those bounds could be used to improve the results in the continuous case. This work is in progress.

**Acknowledgments.** The authors thank Jérôme Le Rousseau and Gunther Uhlmann for their interest in that work and fruitful discussions.

## Appendix

We check in this appendix that Assumptions 2 and 3 hold in several relevant situations.

**A finite-element method** Set $n$ the cardinal of $\mathcal{K}_h$, we shall identify $C(\mathcal{K}_h)$ and $\mathbb{C}^n$ via a so-called numbering of the nodes. We suppose that the discretization of the problem leads to a linear system of the form

$$K_h u_h + q_h u_h = 0, \tag{A.1}$$

where $K_h$ is the $n \times n$ rigidity matrix representing the Laplacian and $q_h u_h$ is the multiplication of the function $q_h$ by the function $u_h$. Note that, for Lagrangian finite elements, the discretization would rather yield a system of the form

$$K_h u_h + M_h[q_h] u_h = 0, \tag{A.2}$$

where $M_h[q_h]$ is the mass matrix, a $n \times n$ matrix whose coefficient depend linearly on $q_h$ defined by

$$(M_h[q_h])_{ij} = \int_\Omega q \phi_i \phi_j,$$

where $(\phi_i)$ is the basis function associated to the node $i$. We suppose that the so-called "mass lumping" technique has been used. This technique reduces the mass matrix to a diagonal matrix, hence to a multiplication coefficient-wise by a vector, which yields a system of the form (A.1).

For each $x, y \in C(\mathring{\mathcal{K}}_h)$, we denote $K_h(x, y)$ the coefficient of the matrix $K_h$ that links the degrees of freedom of $u_h$ associated to the nodes $x$ and $y$. We suppose that $K_h$ obeys the following standard conditions :

- $K_h$ is symmetric, that is $K_h(x, y) = K_h(y, x)$.

- The constants functions of $C(\mathcal{K}_h)$ are in the kernel of $K_h$, that is $K_h(1) = 0$ on $\mathring{\mathcal{K}}_h$.
- $K_h(x, y) = 0$ if $x \not\sim y$ or $x \neq y$

Note that these conditions are natural and usually satisfied when discretizing the Laplace operator.

In order to simplify the notations, define, for any $x \in \mathring{\mathcal{K}}_h$,

$$x^{\pm i/2} = x \pm \frac{h}{2} e_i \quad \text{and} \quad x^{\pm i} = x \pm h e_i$$

The third condition states that

$$K_h u_h(x) = \sum_{y \in \mathcal{K}_h} K_h(x, y) u_h(y)$$

$$= K_h(x, x) u_h(x) + \sum_{i=1}^{k} K_h(x, x^{+i}) u_h(x^{+i}) + K_h(x, x^{-i}) u_h(x^{-i}).$$

The second condition states that

$$K_h(x, x) = -\sum_{i=1}^{k} \left( K_h(x, x^{+i}) + K_h(x, x^{-i}) \right).$$

For each direction of connection $i$, define $\sigma_h^i$ as $\sigma_h^i(x) := h^2 K_h(x^{+i/2}, x^{-i/2})$

Then, for any $u_h \in C(\mathring{\mathcal{K}}_h)$,

$$(\sum_i d_i \sigma^i d_i u)(x) = \frac{1}{h} \sum_{i=1}^{k} (\sigma^i d_i u)(x^{+i/2}) - (\sigma^i d_i u)(x^{-i/2})$$

$$= \frac{1}{h^2} \sum_{i=1}^{k} \left( \sigma_h^i(x^{+i/2})(u_h(x^{+i}) - u_h(x)) - \sigma_h^i(x^{-i/2})(u_h(x) - u_h(x^{-i})) \right)$$

$$= \sum_{i=1}^{k} \left( K_h(x^{+i}, x)(u_h(x^{+i}) - u_h(x)) - K_h(x, x^{-i})(u_h(x) - u_h(x^{-i})) \right)$$

$$= \left( \sum_{i=1}^{k} \left( K_h(x, x^{+i}) u_h(x^{+i}) + K_h(x, x^{-i}) u_h(x^{-i}) \right) \right)$$

$$- u_h(x) \left( \sum_{i=1}^{k} K_h(x, x^{+i}) + K_h(x, x^{-i}) \right)$$

$$= K_h u_h(x).$$

Below, we exhibit some particular cases of interest that verify Assumption 3.

**First instance of Assumption 3.** First, suppose that the mesh $\mathcal{M}_h \in \mathbb{R}^d$ is regular (i.e, $F = Id$) and that the discretization method is the $2d+1$ point Laplacian. In this case, taking $e_1, ..., e_d$ as the canonical vectors of $\mathbb{R}^d$, the matrix $K_h$ is the following:

$$K_h(x, y) = \begin{cases} \frac{1}{h^2} & \text{if } x \sim y \text{ and } x \neq y \\ -\frac{2d}{h^2} & \text{if } x = y \\ 0 & \text{if } x \not\sim y \text{ or } x \neq y. \end{cases}$$

Such a $K_h$ verifies the hypotheses needed in the proof of Theorem 2. Moreover, $\sigma_h^i$ is a constant equal to $\sigma^i = 1$ and $\epsilon_a = \epsilon_d = 0$.

**Second instance of Assumption 3.** If one discretizes the regular mesh with an higher order finite difference Laplacian, say for instance the 9 points Laplacian in two dimensions, it is sufficient to take a higher number of connections. For the 9 points Laplacian, take $e_1 = (1,0)$, $e_2 = (0,1)$, $e_3 = (1,1)$, $e_4 = (1,-1)$, the matrix $K$ takes 4 different values

$$K_h(y,x) = K_h(x,y) = \begin{cases} \frac{a}{h^2} & \text{if } y = x \pm e_1 \text{ or } y = x \pm e_2 \\ \frac{b}{h^2} & \text{if } y = x \pm e_3 \text{ or } y = x \pm e_4 \\ -\frac{4}{h^2}(a+b) & \text{if } x = y \\ 0 & \text{if } x \nsim y \text{ or } x \neq y. \end{cases}$$

By construction of $\sigma$, we have : $\sigma_h^1 = \sigma_h^2 = a$ and $\sigma_h^3 = \sigma_h^4 = b$. Hence $\sigma_h$ has the required properties with $\epsilon_a = \epsilon_d = 0$ for $a + b = 1$

**Third instance of Assumption 3.** We now suppose that the mesh $\mathcal{M}_h \in \mathbb{R}^2$ is a perturbation of $\mathcal{W}_h$, and that the discretrization used is the $P_1$ finite element on triangles. we do not treat the case of tetrahedrons in $\mathbb{R}^3$ which only adds complexity to the notations. It is possible to uniformly refine a cartesian grid into a triangular mesh without adding point by considering the connexions $e_1 = (1,0)$, $e_2 = (0.1)$, $e_3 = (1,1)$ (See Figure 7 left: The black point represent a node $x$, and all the grey points represent the neighbours of $x$).

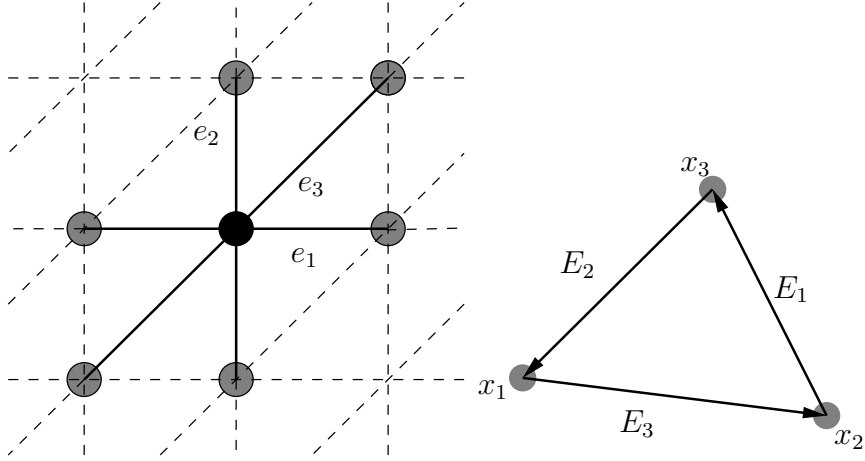

**Figure A1.** Left : The mesh refinement of the rectangular grid for triangular finite elements. Right: Notation on a standard triangle

We recall that, in $P1$ finite elements for triangles, the coefficient $K_h(x,y)$ of the rigidity matrix is defined as

$$K_h(x,y) = -\sum_{T \in T(x,y)} \int_T \nabla \phi_T^x \nabla \phi_T^y,$$

where $T(x,y)$ is the set of triangles that contains both $x$ and $y$ and $\phi_T^x$ is the only linear form of $T$ that is equal to 1 on the vertex $x \in T$ and equal to 0 on the two other vertices of $T$. If the triangle vertices $x_i$ and the triangle edges $E_i$ are numbered as shown in Figure 7 on the right, if $a^\perp$ is the vector $a$ rotated by $\pi/2$, then:

$$\phi^{x_1}(x) = 1 + \frac{E_1^\perp \cdot (x - x_1)}{E_1^\perp \cdot E_2} \quad \phi^{x_2}(x) = 1 + \frac{E_2^\perp \cdot (x - x_2)}{E_1^\perp \cdot E_2}.$$

*Uniform stability estimates for the discrete Calderón problems* 37Consequently:
$$\int_T \nabla \phi_T^{x_1} \nabla \phi_T^{x_2} = \int_T \frac{E_1^\perp \cdot E_2^\perp}{(E_1^\perp \cdot E_2)^2} = 2\frac{E_1 \cdot E_2}{|E_1^\perp \cdot E_2|}$$

We perform an explicit computation of $K(x, x+he_1)$ for which $T(F(x), F(x+he_1))$ is made of two triangles. They both of course have $F(x), F(x+he_1)$ as vertices, the third vertex being either $F(x+he_3)$ or $F(x-he_2)$. For instance, for the $F(x), F(x+he_1), F(x+he_3)$ triangle, we get the following contribution in $\sigma$, that we call $A$ hereafter:

$$A := \int_T \nabla \phi_T^x \nabla \phi_T^{x+he_1}$$
$$= -\frac{1}{2}\frac{(F(x+he_3) - F(x+he_1)) \cdot (F(x+he_3) - F(x))}{|(F(x+he_3) - F(x+he_1))^\perp \cdot (F(x+he_3) - F(x))|}.$$

If we suppose that $F = Id + g$, with $\epsilon_a = \|g\|_{C^1} \leq 1$ a Taylor expansion proves that $A = -1/2$ and summing the different contributions on the triangles, we have

$$\left|\sum \sigma_j e_j \otimes e_j - Id\right| = O(\epsilon_a).$$

Besides, setting $\epsilon_d = \|g\|_{C^2}$, the contribution to $d_i(\sigma^j)$ is of the form
$$d_i A = O(\epsilon_d).$$

For instance, Theorem 2.7 ensures that the scaling that should be compared are $\epsilon_d, \epsilon_a^{1/2}$ and $h^{1/2}$, meaning that if we want $h^{1/2}$ to be the dominant scale, the $C^1$-norm of $g$ has to be smaller than $h$ whereas its $C^2$-norm can scale as $h^{1/2}$ at most. Note that the scaling $\|g\|_{C^2}/\|g\|_{C^1}$ is a scaling of the oscillations of $g$ and that it must be bounded by $h^{-1/2}$.

[1] G. Alessandrini. Stable determination of conductivity by boundary measurements. *Appl. Anal.*, 1988.
[2] Kari Astala and Lassi Päivärinta. Calderón's inverse conductivity problem in the plane. *Ann. of Math. (2)*, 163(1):265–299, 2006.
[3] Kari Astala, Lassi Päivärinta, and Matti Lassas. Calderón's inverse problem for anisotropic conductivity in the plane. *Comm. Partial Differential Equations*, 30(1-3):207–224, 2005.
[4] L. Borcea, V. Druskin, and A. V. Mamonov. Circular resistor networks for electrical impedance tomography with partial boundary measurements. *Inverse Problems*, 26(4), 2010.
[5] Liliana Borcea. Electrical impedance tomography. *Inverse Problems*, 18(6):99–136, 2002.
[6] F. Boyer, F. Hubert, and J. Le Rousseau. Discrete Carleman estimates for elliptic operators and uniform controllability of semi-discretized parabolic equations. *J. Math. Pures Appl. (9)*, 93(3):240–276, 2010.
[7] F. Boyer, F. Hubert, and J. Le Rousseau. Discrete carleman estimates for elliptic operators in arbitrary dimension and applications,. *SIAM J. Control Optim.*, 48:5357–5397, 2010.
[8] F. Boyer, F. Hubert, and J. Le Rousseau. Uniform null-controllability properties for space/time-discretized parabolic equations. *Numer. Math.*, to appear.
[9] Russell M. Brown and Rodolfo H. Torres. Uniqueness in the inverse conductivity problem for conductivities with 3/2 derivatives in $L^p$, $p > 2n$. *J. Fourier Anal. Appl.*, 9(6):563–574, 2003.
[10] A. L. Bukhgeim. Recovering a potential from Cauchy data in the two-dimensional case. *J. Inverse Ill-Posed Probl.*, 16(1):19–33, 2008.
[11] A. P. Calderón. On an inverse boundary value problem. *Seminar on Numerical Analysis and its Applications to Continuum Physics (Rio de Janeiro)*, pages 65–73, 1980.
[12] Edward B. Curtis and James A. Morrow. Determining the resistors in a network. *SIAM J. Appl. Math.*, 50(3):918–930, 1990.
[13] David Dos Santos Ferreira, Carlos E. Kenig, Mikko Salo, and Gunther Uhlmann. Limiting Carleman weights and anisotropic inverse problems. *Invent. Math.*, 178(1):119–171, 2009.
[14] David Dos Santos Ferreira, Carlos E. Kenig, Johannes Sjöstrand, and Gunther Uhlmann. Determining a magnetic Schrödinger operator from partial Cauchy data. *Comm. Math. Phys.*, 271(2):467–488, 2007.